\documentclass{amsart}
\usepackage{amsmath,amsthm,amssymb}

\usepackage{color}
\usepackage{tikz}
\usepackage{tikz-cd}

\usetikzlibrary{cd,positioning}

\usepackage{ytableau}

\usepackage[mathscr]{euscript}

\newcommand{\ad}{\mathrm{ad}}

\newcommand{\af}{{\mathrm{af}}}

\newcommand{\al}{\alpha}

\newcommand{\Aut}{\mathrm{Aut}}

\newcommand{\bGr}{\overline{\Gr}}

\newcommand{\C}{\mathbb{C}}
\newcommand{\cc}{e}
\newcommand{\cI}{\mathscr{I}}

\newcommand{\cl}{\mathrm{cl}}

\newcommand{\cZ}{\mathscr{Z}}

\newcommand{\Diag}{\Gamma_\af}

\newcommand{\eps}{\varepsilon}

\newcommand{\FG}[1]{\pi_1(#1)} 

\newcommand{\fundel}{\sigma} 

\newcommand{\GL}{{GL}}

\newcommand{\Gad}{G_\mathrm{ad}}

\newcommand{\Gm}{\mathbb{G}_m}
\newcommand{\Gr}{\mathrm{Gr}}
\newcommand{\Gsc}{G_\scn}

\newcommand{\hLa}{\hat{\Lambda}}

\newcommand{\Hom}{\mathrm{Hom}}
\newcommand{\hL}[1]{\hat{\Lambda}_{#1}}
 \newcommand{\hLGL}{\hat{\Lambda}_{(n)}}
 \newcommand{\hLSL}{{\hat{\Lambda}_{(n)}}^\circ}

\newcommand{\id}{\mathrm{id}}
\newcommand{\Iaf}{\tilde{I}}

\newcommand{\Is}{\Iaf^s_G}
\newcommand{\isomap}{\stackrel{\cong}{\longrightarrow}}

\newcommand{\KGL}{\mathbb{K}_{\mathrm{GL}_n}}

\newcommand{\K}{\mathbb{K}_{\Gsc}}
\newcommand{\kconj}{\tilde{\omega}}

\newcommand{\Kx}{\mathbb{K}_{G}}
\newcommand{\Kxrat}{\Kx^\Delta}
\newcommand{\ks}[2]{g^{(#1)}_{#2}} 
\newcommand{\tks}[2]{{\tilde{g}}^{(#1)}_{#2}} 
\newcommand{\ksvar}[2]{g^{(#1)}_{#2}(y|b)} 
\newcommand{\tksvar}[2]{{\tilde{g}}^{(#1)}_{#2}(y|b)} 
\newcommand{\la}{\lambda}
\newcommand{\La}{\Lambda}

\newcommand{\LL}{\mathbb{L}_{\Gsc}}
\newcommand{\LGL}{\mathbb{L}_{\GL_n}}

\newcommand{\Lx}{\mathbb{L}_G}
\newcommand{\Lxrat}{\Lx^\Delta}
\newcommand{\map}{\beta}

\renewcommand{\O}{\mathscr{O}}
\newcommand{\om}{\varpi}

\newcommand{\Par}{\mathscr{P}}
\newcommand{\pair}[2]{\langle #1\,,\,#2\rangle}

\newcommand{\PGL}{{PGL}}

\newcommand{\Q}{\mathbb{Q}}

\newcommand{\reg}{\Delta}
\newcommand{\Rep}{{R(T)}}

\newcommand{\RepGL}{{R(\TGL)}}
\newcommand{\Reprat}{\Rep^\Delta}

\newcommand{\rot}{\omega} 

\newcommand{\SL}{{SL}}

\newcommand{\scn}{\mathrm{sc}}
\newcommand{\sh}{{\tau}}

\newcommand{\TGL}{T}
\newcommand{\TSL}{T_{\circ}}
\newcommand{\TPGL}{\overline{T}}

\newcommand{\Tsc}{T_\scn} 

\newcommand{\uu}{v} 
\newcommand{\Waf}{\hat{W}_G}

\makeatletter
\newcommand{\oset}[3][0ex]{%
  \mathrel{\mathop{#3}\limits^{
    \vbox to#1{\kern-2\ex@
    \hbox{$\scriptstyle#2$}\vss}}}}
\makeatother

\newcommand{\Z}{\mathbb{Z}}

\newcommand{\Wafgr}{\Waf^0} 
\newcommand{\WGL}{\tilde{W}_{\GL_n}} 
\newcommand{\WSL}{\tilde{W}_{\SL_n}}
\newcommand{\WGLgr}{\tilde{W}_{\GL_n}^0}
\newcommand{\WSLgr}{\hat{W}_{\SL_n}^0}

\newcommand{\Wxgr}{\Wx^0} 
\newcommand{\Wx}{\tilde{W}_{G}}

\newcommand{\Wfin}{{W_G}}

\newtheorem{thm}{Theorem}[section]
\newtheorem{lem}[thm]{Lemma}
\newtheorem{prop}[thm]{Proposition}
\newtheorem{cor}[thm]{Corollary}

\theoremstyle{remark}
\newtheorem{defn}[thm]{Definition}
\newtheorem{rem}[thm]{Remark}
\newtheorem{ex}[thm]{Example}

\usepackage{chngcntr}

\numberwithin{equation}{section}

\title[Equivariant $K$-homology of affine Grassmannian]{Equivariant $K$-homology of affine Grassmannian and $K$-theoretic double $k$-Schur functions}

\date{}

\author{Takeshi Ikeda}
\address{Department of Pure and Applied Mathematics, School of Fundamental Science and Engineering, Waseda University, 3-4-1 Okubo, Shinjuku-ku, Tokyo 169-8555, Japan}
\email{gakuikeda@waseda.jp}
\author{Mark Shimozono}
\address{Department of Mathematics, 460 McBryde Hall, Virginia Tech, 255 Stanger St., Blacksburg, VA,
24601, USA}
\email{mshimo@math.vt.edu}
\author{Kohei Yamaguchi}
\address{Graduate School of Mathematics, Nagoya University, Furo-cho, Chikusa-ku, Nagoya 464-8602, Japan}
\email{yamaguchi.kohei@math.nagoya-u.ac.jp}

\begin{document}
\begin{abstract}
We study the torus equivariant $K$-homology ring of the affine Grassmannian $\Gr_G$ where $G$ is a connected reductive linear algebraic group. In type $A$, we introduce equivariantly deformed symmetric functions called the $K$-theoretic double $k$-Schur functions as the Schubert bases. The functions are constructed by Demazure operators acting on equivariant parameters. As an application, we provide a Ginzburg-Peterson type realization of the torus-equivariant $K$-homology ring of $\Gr_{\SL_n}$ as the coordinate ring of a centralizer family for $\PGL_n(\C)$. 
\end{abstract}
\maketitle

\tableofcontents

\section{Introduction}

Let $G$ be a simple simply connected complex algebraic group, $T$ a maximal torus of $G$.
Let $\Gr_G=G(\mathbb{F})/G(\mathbb{O})$ be the affine Grassmannian of $G$, where $\mathbb{O}=\C[\![t]\!],\, \mathbb{F}=\C(\!(t)\!)$.
 The $T$-equivariant $K$-homology $K_*^T(\Gr_G)$ 
of $\Gr_G$ 
has the structure of a commutative and cocommutative $\Rep$-Hopf algebra \cite{LSS:K}, where $\Rep$ is the representation ring of $T$.
It has two distinguished $\Rep$-bases $\{\O_w\}$ and $\{\cI_w\}$ consisting of the respective classes of the Schubert structure sheaves and the ideal sheaves of the boundary of the Schubert variety. 
These bases are indexed by the set of the affine Grassmann elements, i.e. the minimum-length coset representatives in the affine Weyl group of $G$ with respect to the classical Weyl group.  


We are interested in the \emph{Schubert calculus} for $K_*^T(\Gr_G)$, i.e. the
explicit description of the structure constants for the basis $\{\O_w\}$.
The Schubert calculus for the (equivariant) homology ring $H_*^T(\Gr_G)$ has been studied intesively (see \cite{LLMSSZ} and references therein),
since Peterson \cite{Pet} (see also \cite{LS:Q=Aff}) introduced the ``quantum equals affine'' principle, which asserts that,  after harmless localization, the equivariant quantum cohomology ring $QH_T^*(G/B)$ of the flag variety $G/B$ is isomorphic to the equivariant homology ring $H^T_*(\Gr_G)$, under which 
Schubert classes correspond, establishing the equality of structure constants on both sides.
The $K$-theoretic analogue of Peterson's isomorphism in type $A$ non-equivariant case was studied in
\cite{IIM} (see also \cite[\S 4]{IIN}). 
The general type correspondence was 
conjectured in \cite{LLMS} and proved by Kato \cite{Kat}. 
Chow and Leung \cite{CL} gave another proof. 


The main object of study in this paper is $K_*^T(\Gr_{G})$ for type $A$ groups $G$ including $\GL_n(\C)$. 
In the course of study we generalize the theory of the \emph{affine level zero $K$-nil-Hecke algebra} $\Kx$, developed in \cite{LSS:K}, 
to an arbitrary reductive algebraic group $G$. 
Such a generalization, in particular working with a nontrivial fundamental group $\FG{G}$, also has interesting features in its own right. 
One such benefit is a clear explanation of 
the $k$-rectangle factorization property (see \cite{LM} for the analogous property of $H_*(\Gr_{\SL_n})$ and the $k$-Schur basis) of the structure sheaves $\mathscr{O}_w$ and its generalization to reductive $G$ (\S \ref{SS:k rectangles}). We give a realization of $K_*^T(\Gr_{\GL_n})$
as an $\Rep$-algebra of certain 
equivariantly deformed symmetric functions, and 
introduce 
the \emph{$K$-theoretic double $k$-Schur functions} as the
algebraic counterparts of the Schubert bases for both $\{\O_w\}$ and $\{\cI_w\}$. 
As an application we provide another algebraic model
of $K_*^T(\Gr_G)$ for type $A$ groups 
as a coordinate ring of certain centralizer family $\mathscr{Z}_{G^\vee}$
of the dual group $G^\vee$.

\subsection{$K$-theoretic affine level zero nil-Hecke algebra for reductive $G$}
Let $G$ be a reductive complex algebraic group, $T$ a maximal torus of $G$, and $B$ a Borel subgroup of $G$ containing $T$.
Let $\Wx=X^\vee\rtimes \Wfin$ be the extended 
 affine Weyl group, where $X^\vee$ is the lattice of cocharacters of $T$, and $W_G$ the Weyl group of $G.$
The group is isomorphic to $\FG{G}\ltimes \hat{W}_G$, where 
$\hat{W}_G$ is the (non-extended) affine Weyl group generated by the standard Coxeter generators $s_i\;(i\in \Iaf)$, where 
$I$ and $\Iaf=I\cup \{0\}$ are the index sets of the Dynkin diagram of $G$ and its untwisted affine extension. $\Wx$ acts on $\Rep$ by the level zero action: $X^\vee$ acts by the identity and $\Wfin$ acts naturally. In particular $s_0$ acts by $s_\theta$ where $\theta$ is the highest root.
We denote by $\alpha_i$ the simple root corresponding to $i\in I$. We also set $\alpha_0=-\theta.$
For $i\in \Iaf$ we define the Demazure operator 
$$
T_i=({1-e^{\alpha_i}})^{-1}(s_i-1)$$
as an element of a twisted group ring of $\tilde{W}_G$ with coefficients in a suitable localization of $\Rep$; see \S \ref{SS: def Kx}.
We also use $D_i:=T_i+1.$
It is known that the $T_i$ and $D_i$ satisfy the braid relations, so that $T_w$ and $D_w$
are defined for $w\in \hat{W}_G$. 
We also define $T_w = \fundel T_v$ and $D_w=\fundel  D_v$ for $w=\sigma v\,(\sigma\in \FG{G}, \, v\in \hat{W}_G).$
Let $\Kx$ the left $\Rep$-module generated by $T_w$ for $w\in \tilde{W}_G$.
In fact, $\Kx$ is a ring with $\Rep$-basis $\{T_w\mid w\in \Wx\}$
called the \emph{affine level zero $K$-nil-Hecke ring}.

The ring $\Kx$ has a commutative subring $\Lx$, the \emph{$K$-theoretic Peterson subalgebra}, defined as the centralizer 
of $\Rep$ in 
$\Kx$. 
There is an $\Rep$-Hopf algebra isomorphism
$K_*^T(\Gr_G)\cong \Lx$ \cite{LSS:K}. Denote by $k_w$ and $\ell_w$ be the images in $\Lx$ of the ideal sheaf $\cI_w$ and structure sheaf $\O_w$ respectively. The bases $\{k_w\}$ and $\{\ell_w\}$
have algebraic characterizations \cite{LSS:K}, \cite{LLMS}.

One of the key constructions of this paper is an action $*$ of $\Kx$ on $\Lx$, whose homological analogue is due to Peterson \cite{Pet}.
Let $\Wxgr$ be the set of affine Grassmann elements in $\Wx$, the minimum-length coset representatives in $\Wx/\Wfin$.  The $*$-action allows the following simple
definition of the algebraic Schubert bases $\{k_w\}$ and $\{ \ell_w\}$ of $\Lx$ for $w\in \tilde{W}^0_G$ (see Prop. \ref{P:star basis formula}):
\begin{equation}
\label{E: k and l}
    k_w = T_w * 1, \quad
    \ell_w = D_w * 1.
\end{equation}
This construction should play a fundamental role in the study of the Schubert calculus of $K^T_*(\Gr_G)$.
Its homological analogue appears in \cite{Pet} but the $*$-action is not explicitly given in \cite{LSS:K}.

\subsection{$K^T_*(\Gr_G)$ as equivariantly deformed ring of symmetric functions}

In \cite{LSS:K} a family of symmetric functions 
called \emph{$K$-theoretic $k$-Schur functions},
was introduced as a counterpart of the non-equivariant Schubert 
basis 
of $K_*(\Gr_{\SL_n})$.
More precisely, the function $g_w^{(k)}(y)$ defined in \cite{LSS:K}
corresponds to the non-equivariant ideal sheaf, which we denote by  $\cI_w^0$. 
The function $\tilde{g}_w^{(k)}(y)$, called the \emph{closed $K$-theoretic $k$-Schur function},
corresponding to $\O_w^0$ 
was studied in \cite{Tak}, \cite{BMS}, \cite{IIN}.  
The construction of $g_w^{(k)}(y)$ in \cite{LSS:K} used a non-equivariant version of $\mathbb{L}_{G}$ called the \emph{$K$-theoretic affine Stanley algebra}. 
In this article we work directly in the equivariant setting, using divided difference operators to create the Schubert bases in symmetric functions. This approach was inspired by Nakagawa and Naruse \cite{NN} and applied \cite{LLS:backstable} in the infinite rank type $A$ case.
For $G$ of type $A$, we define an injective $R(T)$-Hopf morphism from $\Lx$ to an equivariantly deformed algebra of symmetric functions, so that this map is $\Kx$-equivariant.  We transport the bases $k_w$ and $\ell_w$ to symmetric functions and define their image to be the
\emph{$K$-theoretic double $k$-Schur functions}.

To explain the construction in more detail, let us fix some notation. 
Let $G=\GL_n(\C)$ and $T\subset G$ the maximal torus of diagonal matrices.
The cocharacter lattice is $X^\vee=\oplus_{i=1}^n\Z \eps_i$ (see Example \ref{X:GL root data}).
We have the fundamental group $\FG{G}\cong \Z$, for which we fix a generator $\sh$ (see Example \ref{X: extended affine Weyl type A}).
We have
$\tilde{W}_{\GL_n}\cong \FG{G} \ltimes \tilde{W}_{\SL_n}$.

For any commutative ring $R$, let $\hat{\Lambda}_{R}$ be the completed ring 
of symmetric functions in 
the variables $y=(y_1,y_2,\ldots)$
with coefficients in $R$; see \S \ref{S: Def K-k-Schur}.
The representation ring of $T$ is $\Rep=\Z[e^{\pm a_1},\ldots,e^{\pm a_n}],$ where $\{a_i\}$
is the $\Z$-basis of the lattice $X$ of characters of $T$ that is dual to the basis $\{\eps_i\}$ of $X^\vee$ under the canonical pairing. Let $b_i=1-e^{-a_i}\in \Rep$ for $1\le i\le n.$
We define a left $\mathbb{K}_{\GL_n}$-module structure on $\hL{\Rep}$ in Proposition \ref{P: KGL acts on hat Lambda}.
Under this action the symmetric group $W_{\GL_n}$ acts by permuting the variables $b_1,\ldots,b_n$, and the translation element $t_{\eps_i}\in\tilde{W}_{\GL_n}$ acts on $\hL{\Rep}$
by multiplication by the symmetric series
\begin{equation}
\label{def:Omega}
\Omega(b_i| y):=\prod_{j=1}^\infty(1-b_i y_j)^{-1}.
\end{equation}
For $w\in \WGLgr$, define 
\begin{equation}
\tksvar{k}{w}=D_w(1),\quad 
\ksvar{k}{w} =T_w(1). 
\end{equation}
We call $\ksvar{k}{w}$ the 
\emph{
$K$-theoretic double $k$-Schur function} and
$\tksvar{k}{w}$ the \emph{closed
$K$-theoretic double $k$-Schur function}.

Let $\hLGL$ be the $\Rep$-span of $\tksvar{k}{w}\;(w\in \WGLgr)$.
There is a natural $\Rep$-Hopf algebra structure on $\hLGL$ (\S \ref{S: Def K-k-Schur}).

\begin{thm}\label{thm:main}
There is a $\KGL$-module and $\RepGL$-Hopf isomorphism 
\begin{align}
 \label{E:the iso}
&K_*^{\TGL}(\Gr_{\GL_n})\isomap \hLGL ,\\
\O_w \mapsto &\tksvar{k}{w}, \quad
\cI_w \mapsto g^{(k)}_w(y|b)\quad \text{for $w\in \WGLgr$.}
\end{align}
\end{thm}

See \S \ref{SS:SL and PGL}, for the corresponding statement for $G=\SL_n(\C)$ and $\PGL_n(\C).$
In particular, the $K$-theoretic double $K$-Schur function $g_w^{(k)}(y|b)$ specializes to the $K$-theoretic 
$k$-Schur function $g_{w}^{(k)}(y)$ at $b=0$ (Corollary \ref{C:forget to K k schur}). 

\subsection{$K^T_*(\Gr_{G})$ as the coordinate ring of a centralizer family}
A theorem of Ginzburg \cite{Gi} and Peterson \cite{Pet} (see also \cite{YZ}) 
asserts that the equivariant homology ring $H^T_*(\Gr_G)$, is isomorphic to the coordinate ring of a centralizer family in the Lie algebra $\mathfrak{g}^\vee$ that is the Langlands dual to the Lie algebra of $G$. In type $A$ this was worked out explicitly \cite{LS:Toda}, \cite{LS:double Toda}, as well as an explicit substitution producing $k$-Schur functions from quantum Schubert polynomials.
We give a similar centralizer construction
in $K$-homology for $G=\SL_n(\C)$ and $G^\vee=\PGL_n(\C).$
Consider the following matrix with entries on $\Rep$
\begin{align}\label{E:matrix A}
 A = \begin{pmatrix}
  e^{-a_1} & -1& &  & \\
   & e^{-a_2} &-1 &&\\
&& \ddots&\ddots&\\
&& & \ddots&-1\\
&&&&e^{-a_n}
\end{pmatrix}.
\end{align}
Since $e^{-a_i}$ is a function on the maximal torus $T$ of $\GL_n(\C),$
this gives a family of elements in $\GL_n(\C)$ parametrized by 
$T$.
Let $\cZ_{\GL_n}$ be the closed subscheme of $T\times \GL_n(\C)$
consisting of the pair $(A,Z)\in T\times \GL_n(\C)$ such that $AZ=ZA.$
For such a pair $Z$ must be upper-triangular.
$\cZ_{\GL_n}$ is an affine scheme over $T$ by the first projection.
The coordinate ring $\O(\cZ_{\GL_n})$ has a natural structure of an $R(T)$-Hopf
algebra. Let 
$z_{ij}\,(1\le i\le j\le n)$ be the regular function on $\cZ_{\GL_n}$ corresponding the $(i,j)$ entry of 
$Z$. Let $\omega$ be the automorphism of $\Rep$ such that $\omega(b_i)=b_{i+1}$ for $1\le i\le n$ where $b_{n+1}=b_1$.

\begin{thm} There is an $\Rep$-Hopf isomorphism
$\O(\cZ_{\GL_n}) \cong \hLGL$ sending
\begin{align}
z_{ii}&\mapsto\Omega(b_i|y)&\qquad&\text{for $1\le i\le n$,} \\
z_{ij}&\mapsto e^{a_i+\cdots+a_{j-1}}\Omega(b_i|y)
g_{\rho_{j-i}}^{(k)}(y|\omega^i(b)) &&\text{for $1\le i<j\le n$,}
\end{align}
where $\rho_l=s_{l-1}\cdots s_1s_0$ for $1\le l\le n-1$.
\end{thm}

Let $\cZ_{\GL_n}|_{\TSL}$ denote the restriction to $\TSL\subset T.$
We define $\cZ_{\PGL_n}=(\cZ_{\GL_n}|_{\TSL})/\Gm$, where $\Gm=\C^\times$
acts on $\cZ_{\GL_n}$ by dilation of $Z.$ 
$\O(\cZ_{\PGL_n})$ is generated by $z_{ij}/z_{11}$ as an $R(\TSL)$-algebra.
\begin{cor} There is an $R(\TSL)$-Hopf isomorphism
$\O(\cZ_{\PGL_n})\cong K_*^{\TSL}(\Gr_{\SL_n})$.
\end{cor}

\subsection*{Related work and further questions}
Based on the type $A$ result in this paper, the ``quantum equals affine'' perspective in type $A$ will be studied in \cite{IINY}. The integrable system called the
\emph{relativistic Toda lattice} plays the role of giving 
a bridge between the affine and quantum sides. 
This work is a natural extension of  \cite{LS:Toda}, \cite{LS:double Toda}  
for (co)homology, and also of \cite{IIM} on the non-equivariant quantum $K$-theory.

  It will be interesting to study symmetric function realizations of equivariant bases for the homology of affine Grassmannians of other groups. For type $C$, the non-equivariant (co)homology case was studied by Lam, Schilling, and Shimozono \cite{LSS:C}.
The extensions of this result to equivariant homology and equivariant $K$-homology 
will be studied in forthcoming papers \cite{IIS} and \cite{IISY}.

We expect that there is a similar centralizer construction of $K_*^T(\Gr_G)$ 
and symmetric function constructions for general reductive groups $G$.
The $G(\mathbb{O})$-equivariant
 $K$-homology ring
$K_*^{G(\mathbb{O})}(\Gr_G)$ was 
studied by Bezrukavnikov, Finkelberg, and Mirkovi\'c \cite{BFM}.
They identify $K_*^{G(\mathbb{O})}(\Gr_G)$
with the coordinate ring of the universal centrlizer.

It should be an interesting problem to give an equivariant version of the 
$K$-theoretic Catalan formula for $g_w^{(k)}(y)$ 
due to Blasiak, Morse, and Seelinger \cite{BMS}, and for $\tilde{g}_w^{(k)}(y)=\tilde{g}^{(k)}_w(y|0)$ by \cite{IIN}.

\subsection*{Organization} In Section \ref{S:K-nil}, we discuss
the affine level zero $K$-nil-Hecke ring for 
arbitrary reductive group $G.$
The fundamental result is  
Theorem \ref{T:G iso} establishing 
an isomorpshim between 
$K_*^T(\Gr_G)$ and the Peterson subalgebra $\Lx.$
The $*$ action of $\Kx$ on $\Lx$ is given in \S \ref{SS:star action}.
In \ref{SS:k rectangles} we discuss the factorization property. 
In \S \ref{SS:adjoint} we study the case that $G$
is of adjoint type. 
In Section \ref{S:K k Schur}, we introduce the $K$-theoretic double $k$-Schur functions. The $k$-rectangle factorization is
given in \S \ref{SS: factorization}.
In Section \ref{S:cent}, we give the centralizer constructions of $K_*^T(\Gr_G)$ for type $A$
groups. In Appendix \ref{S:Molev}, we discuss the 
infinite rank analogue of the $K$-theoretic doule $k$-Schur functions.
In Section \ref{S:double k schur}, we explain the 
double $k$-Schur functions and its infinite rank analogue 
called Molev's dual Schur functions.
\subsection*{Acknowledgements}
We thank Noriyuki Abe, Shinsuke Iwao, Hajime Kaji, \newline Thomas~Lam, and Satoshi Naito for helpful discussions.
This work was supported by JSPS KAKENHI Grant Numbers 23H01075, 	22K03239, 	20K03571,	20H00119.


\section{$K$-theoretic affine nil-Hecke algebra for reductive $G$}\label{S:K-nil}
In this section we generalize the fundamental construction
in \cite{LSS:K} to reductive group $G$. We introduce 
the level zero affine $K$-nil-Hecke algabra $\Kx$.
The $T$-equivariant $K$-homology ring $K_*^T(\Gr_G)$ is shown to isomorphic to
a commutative
algebra $\Lx$ of $\Kx$ called the Peterson subalgebra. 

We also introduce a left  
$\Kx$-module structure on $\Lx$ called the star action.
This enable us to construct the Schubert bases of $\Kx$ algebraically. 
In \S \ref{SS:k rectangles}, we prove
a factorization property of $\O_w$, called $k$-rectangle factorization in type $A$,
 with respect to $\O_{t_{-\om_i^\vee}}$.
 \S \ref{SS:adjoint} is devoted to the study of the group of adjointe type. In \S \ref{ssec:Affine Dynkin auto} we discuss symmetries of $K_*^T(\Gr_G)$
 arising from automorphisms of affine Dynkin diagram fixing the zero-th node.
\subsection{Extended affine Weyl group}

Let us denote by $\Gsc$ and $\Gad$
the semisimple algebraic groups
having the same Lie type as $G$ but are of simply-connected and adjoint type respectively.
\subsubsection{Root data and the classical Weyl group $\Wfin$}
Let $X(T)=\Hom(T,\Gm)$ and $X^\vee(T)=\Hom(\Gm,T)$ be the
lattice of characters and cocharacters of $T$ respectively,
where $\Hom$ is understood to be the set of homomorphisms
as algebraic groups, and 
$\Gm=\C^\times$ is the multiplicative group. 
There is a perfect pairing 
\begin{align}\label{E:perfect pairing}
\pair{\cdot}{\cdot}: X^\vee(T)\times X(T) \to \Hom(\Gm,\Gm)\cong\Z
\end{align}
given by $(f,g)\mapsto g\circ f$. 
Let $\Phi\subset X(T)$ and $\Phi^\vee\subset X^\vee(T)$ be the root system and the coroot system of $(G,T)$ respectively.
The quadruple $(X(T),\Phi,X^\vee(T),\Phi^\vee)$
is called the \emph{root datum} of $(G,T)$ (cf. \cite[\S 7.4]{Sp}).
We fix a set of simple roots $\{\alpha_i\mid i\in I\}\subset \Phi$
and a set of simple coroots $\{\alpha_i^\vee\mid i\in I\}\subset \Phi^\vee$,
both indexed by a finite set $I$, 
so that these are positive with respect to $B$.

The root lattice is $Q=\bigoplus_{i\in I} \Z\alpha_i\subset X(T)$ and the coroot lattice is $Q^\vee=\bigoplus_{i\in I} \Z\alpha_i^\vee\subset X^\vee(T)$.
The weight lattice $P\subset \Q\otimes_\Z X(T)$ is the lattice dual to $Q^\vee$
and the coroot lattice $P^\vee\subset \Q\otimes_\Z X^\vee(T)$ is the lattice dual to $Q$.
The $\Z$-basis $\{\om_i^\vee\mid i\in I\}$ of $P^\vee$ dual to the simple roots, consists of the fundamental coweights. The $\Z$-basis $\{\om_i\mid i\in I\}$ of $P$ dual to the simple coroots, consists of the fundamental weights.

The classical Weyl group $\Wfin=N(T)/T$ acts on $T$
and therefore on $X(T)$, $X^\vee(T)$, $Q$, $Q^\vee$, $P$, and $P^\vee$.
The group $\Wfin$ is isomorphic to the subgroup of $\GL(X)$ (resp. of $\GL(X^\vee)$)
generated by the simple reflections $s_i$ for $i\in I$, which act on $X$ (resp $X^\vee$) by
\begin{align}
s_i(\la) &= \la - \pair{\alpha_i^\vee}{\la} \alpha_i \qquad\text{for $\la\in X(T)$} \\
\label{E:s on Xv}
  s_i(\mu) &= \mu - \pair{\mu}{\alpha_i} \alpha_i^\vee \qquad\text{for $\mu\in X^\vee(T)$.}
\end{align}
We have
\begin{align*}
    \pair{w\mu}{w\la} = \pair{\mu}{\la} \qquad\text{for all $w\in \Wfin$, $\la\in X(T)$ and $\mu\in X^\vee(T)$.}
\end{align*}

\begin{ex}\label{X:GL root data}
For $G=\GL_n(\C)$ let $T$ be a maximal torus of $G$. We have 
$$
X(\TGL)=\bigoplus_{i=1}^n\Z a_i,\quad
X^\vee(\TGL)=\bigoplus_{i=1}^n\Z \eps_i
$$
with the pairing $X^\vee(\TGL)\times X(\TGL)\rightarrow \Z$ given by 
$\langle \eps_i,a_j\rangle=\delta_{ij}$.
The root system of 
$(G,T)$ is $$\Phi=\{\pm(a_i-a_j)\mid 1\le i <j\le n\},\;\Phi^\vee=\{\pm(\eps_i-\eps_j)\mid 1\le i < j\le n\}.$$
The classical simple roots (resp. simple coroots) are $\alpha_i = a_i - a_{i+1}$ (resp. $\alpha_i^\vee= \eps_i -\eps_{i+1}$) for $i\in I=\{1,2,\dotsc,n-1\}$. We also use $\alpha_0 = a_n - a_1$.
Note that the root datum $(X(\TGL),\Phi,X^\vee(\TGL),\Phi^\vee)$ is self-dual
via the isomorphism $X(\TGL)\cong X^\vee(\TGL)$ of abelian groups given by 
$a_i\mapsto \eps_i \;(1\le i\le n).$

$\Wfin$ is the symmetric group $S_n$ on $n$ symbols, 
acting on $X^\vee(T)$ and $X(T)$ by permuting the
indices of the $\eps_i$ and $a_i$.

For other groups of type $A$ such as the simply-connected group $\SL_n(\C)$ and 
adjoint group $\PGL_n(\C)$, the simple roots, simple coroots, 
$\Phi$, and $\Wfin$ are all ``the same" but live in slightly different 
lattices; see \S \ref{SS:type A root data} 
for the root data for $\SL_n(\C)$ and $\PGL_n(\C)$.
\end{ex}

\subsubsection{Extended affine Weyl group via translations}
\label{SS:extended affine translations}
In the following, we often denote $X(T)$ (resp. $X^\vee(T)$) by $X$ (resp. $X^\vee$) when 
the chosen maximal torus $T$ is clear from the context.

The \emph{extended affine Weyl group} of $G$ is by definition
\begin{align}\label{E:general affine Weyl}
   \Wx = X^\vee \rtimes \Wfin.
\end{align}
Let $t_\la$ denote the image in $\Wx$ of $\la\in X^\vee$.
We have 
\begin{align}\label{E:translation finite commutation}
w t_\la w^{-1} = t_{w(\la)}\qquad\text{for $w\in \Wfin$ and $\la\in X^\vee$.}
\end{align}
$\Wx$ has a faithful 
realization as a subgroup of affine linear transformations on $X^\vee$:
$\Wfin$ acts on $X^\vee$ via \eqref{E:s on Xv} while
$X^\vee$ acts on itself by translations.

There is a group homomorphism $\cl: \Wx\to \Wfin$ defined by
$\cl(t_\la u)=u$ for $\la\in X^\vee$ and $u \in \Wfin$.
The \emph{level zero action} of $\Wx$ on $X$ is the nonfaithful action defined by
$w(\la)=\cl(w)(\la)$ for $w\in \Wx$ and $\la\in X$. 

The \emph{affine Weyl group} is the subgroup of $\Wx$ defined by
\begin{align*}
   \Waf  = Q^\vee \rtimes \Wfin.
\end{align*}
This is a Coxeter group $\langle s_i\mid i\in\Iaf\rangle$ of affine type for $\Iaf=I\cup\{0\}$. It is generated by $\Wfin=\langle s_i\mid i\in I\rangle$ together with the 
affine simple reflection
\begin{align*}
    s_0 = t_{\theta^\vee} s_\theta
\end{align*}
where $s_\theta\in \Wfin$ and $\theta^\vee\in \Phi^\vee$ are the associated reflection and associated coroot of the highest root $\theta$.

\begin{ex} \label{X:GL translation elements}
For $G=GL_n(\C)$, $\WGL$ is generated by the translation elements,
which are indexed by $X^\vee=\bigoplus_{i=1}^n\Z \eps_i$ (see Example \ref{X:GL root data}),
and the symmetric group $\Wfin=S_n$. We have $\theta^\vee = \eps_1 -\eps_n$,
$\theta = a_1 - a_n$, and $s_\theta$ is the transposition that exchanges $1$ and $n$. The affine Weyl group $\Waf$ is the group associated with the simply-connected group $\SL_n$; it has translation elements from $Q^\vee \subset X^\vee$.
\end{ex}

\subsubsection{Fundamental group and affine Dynkin automorphisms}
\label{SS:fundamental group}
The fundamental group $\FG{G}$ of $G$ acts on the affine Weyl group $\Waf$ by affine Dynkin automorphisms. 
We have a commutative diagram of groups and group homomorphisms, whose explanation follows.

\begin{equation}\label{E:fundamental as Dynkin auto}
\begin{tikzcd}
  X^\vee/Q^\vee \arrow[r,"\overline{\phi}"] \arrow[d,swap,"\cong"]& P^\vee/Q^\vee\arrow[d,"\cong"] & \\
     \FG{G} \arrow[r,dotted] & \Aut(\Diag)^s \arrow[r,"\subset"] & \Aut(\Diag) \arrow[r] & \mathrm{Aut}(\Waf)
\end{tikzcd}
\end{equation}
The restriction of $\mu\in X^\vee=\Hom(\Gm,T)$ to $S^1\subset \Gm$ 
is a continuous map $S^1\to T\hookrightarrow G$, yielding an element of $\FG{G}$. 
One may show this defines a surjective group homomorphism $X^\vee\to\FG{G}$ with kernel $Q^\vee$, inducing an isomorphism of groups \cite[Prop. 1.11]{Bo}
\begin{align}\label{E:fund Xvee/Qvee}
X^\vee/Q^\vee \cong \FG{G}.
\end{align}

Consider the composite map $\phi: X^\vee\to P^\vee$ defined by 
\begin{align}\label{E:Xvee to Pvee}
\begin{tikzcd}[ampersand replacement=\&]
X^\vee \arrow[r,"\cong"] \& \Hom_\Z(X,\Z) \arrow[r,"\text{restrict}"] \& \Hom_\Z(Q,\Z) \arrow[r,"\cong"] \& P^\vee
\end{tikzcd}
\end{align}
where the first isomorphism comes from the perfect pairing \eqref{E:perfect pairing}, the second map is restriction from $X$ to $Q$, and the last map is by definition. Since $\phi$ restricts to the inclusion of $Q^\vee$ into $P^\vee$, there is an induced group homomorphism
$\overline{\phi}: X^\vee/Q^\vee \to P^\vee/Q^\vee$. 

\begin{rem}\label{R:Xvee mod Qvee embeds}
The map $\overline{\phi}$ is injective if $G$ is semisimple and an isomorphism if $G$ is of adjoint type.
\end{rem}

Let $\Aut(\Diag)$ be the automorphism group of the affine Dynkin diagram $\Diag$, which is a subgroup of the group of permutations of $\Iaf$. $\Aut(\Diag)$ acts on $\Waf$ by
\begin{align}\label{E:pi on simple reflections}
  \fundel(s_i) = s_{\fundel(i)} \qquad\text{for $\fundel\in\Aut(\Diag)$, $i\in \Iaf$.}
\end{align}
$\Aut(\Diag)$ similarly acts on the affine root lattice $Q_\af$ by permuting the affine simple roots. 
This action fixes the null root $\delta$. This induces an action of $\Aut(\Diag)$ on $Q\cong Q_\af/\Z\delta$ defined by
\begin{align}\label{E:pi on simple roots}
    \fundel(\alpha_i) = \alpha_{\fundel(i)} \qquad\text{for $i\in \Iaf$ and $\fundel\in\Aut(\Diag)$.}
\end{align}

Let $\Iaf^s := \Aut(\Diag)\cdot\{0\}\subset \Iaf$ be the set of special (cominuscule) nodes. There is a bijection $\Iaf^s\cong P^\vee/Q^\vee$ given by $i\mapsto \om_i^\vee+Q^\vee$ for $i\in \Iaf^s$, where $\om_0^\vee=0$ by convention. Via this bijection, for $i\in \Iaf^s$, 
let $\pi_i$ be the permutation of $\Iaf^s$ induced by the permutation of $P^\vee/Q^\vee$ given by addition by $\om_i^\vee+Q^\vee$. 
This permutation of $\Iaf^s$ extends uniquely to an element of $\Aut(\Diag)^s$ which we also denote by $\pi_i$.
Let $\Aut(\Diag)^s = \{\pi_i\mid i\in \Iaf^s\} \le \Aut(\Diag)$; we call it the group of special affine Dynkin automorphisms. Thus we have a group isomorphism
\begin{align*}
    P^\vee/Q^\vee \cong \Aut(\Diag)^s.
\end{align*}
This group is isomorphic to $\pi_1(\Gad).$
Let $\Is$ be the set of elements $i\in \Iaf^s$ corresponding to elements $\pi_i$ in the image of 
$\overline{\phi}: X^\vee(T)/Q^\vee\to P^\vee/Q^\vee$.                
\begin{align*}\FG{G}\cong X^\vee/Q^\vee\to \{\pi_i\mid \text{$i\in \Is$}\}\le\Aut(\Diag)^s.
\end{align*}

\begin{rem}\label{R:inverse in fundamental group}
Negation in $X^\vee/Q^\vee$ induces the automorphism of the abelian group $\FG{G}$ given by taking inverses.
For $i\in \Is$ define $i^*\in \Is$ by $\pi_i^{-1}=\pi_{i^*}$. 
The map $i\mapsto i^*$ extends to an element of $\Aut(\Diag)$. Explicitly $0^*=0$ and
for $i\in I$ we have $w_0(\om_i^\vee) = - \om_{i^*}^\vee$, where $w_0\in\Wfin$ is the longest element.
\end{rem}

\subsubsection{Extended affine Weyl group via fundamental groups}
Using the homomorphism $\FG{G}\to \mathrm{Aut}(\Waf)$ of \eqref{E:fundamental as Dynkin auto},
$\Wx$ is the semidirect product
\begin{align}\label{E:extends affine weyl}
  \Wx &\cong \FG{G} \ltimes \Waf \\
  \label{E:fundamental affine commutation}
  \fundel w \fundel^{-1} &= \fundel(w)\qquad\text{for $\fundel\in\FG{G}$ and $w\in \Waf$}
\end{align}
where the action of $\fundel$ on $\Waf$ is defined by 
\eqref{E:pi on simple reflections}.

For the expression of $\fundel\in\FG{G}$ in the presentation \eqref{E:general affine Weyl} of $\Wx$, we define the notation
\begin{align}\label{E:Pi as translation finite}
\fundel = t_{\fundel} u_\fundel\qquad\text{with $t_\fundel\in X^\vee$ and $u_\fundel=\cl(\fundel)\in \Wfin$.}
\end{align}
Since $\cl$ is a homomorphism, its restriction to $\FG{G}\le \Wx$ yields a group homomorphism
$\pi_1(G)\to \Wfin$. 
For $G$ semisimple \eqref{E:Pi as translation finite} reads
\begin{align}\label{E:fundamental into translation finite}
    \pi_i = t_{\om_i^\vee} u_i\qquad\text{for $i\in \Is$}
\end{align}
where $u_i$ is the shortest element of $\Wfin$ such that $u_i(\om_{i^*}^\vee)=-\om_i^\vee$;
see Remark \ref{R:inverse in fundamental group}.
We note that 
\begin{equation}\label{E:fundamental inverse into translation finite}
    \pi_i^{-1} = t_{\om_{i^*}^\vee} u_i^{-1}\qquad\text{for $i\in \Is$}.
\end{equation}
\begin{rem} \label{R:extended affine Weyl}
\begin{enumerate}
\item $\Wx$ is the usual extended affine Weyl group 
if and only if $X^\vee=P^\vee$, that is, $G$ is semisimple of adjoint type.
In this case $\Is=\Iaf^s$ and $\FG{G} \cong P^\vee/Q^\vee\cong\Aut(\Diag)^s$.
\item $\Wx=\Waf$ if and only if $X^\vee=Q^\vee$, that is,
$G$ is semisimple and simply-connected. In this case $\Is= \{0\}$ and $\FG{G}=\{\id\}$.
\item If $G$ is semisimple then $Q^\vee\subset X^\vee \subset P^\vee$.
\end{enumerate}
\end{rem}

\begin{ex}\label{X: extended affine Weyl type A}
Typical Type $A$ cases are as follows.
\begin{enumerate} 
\item Since $G=\SL_n(\C)$ is simply-connected,
$\Is = \{0\}$, $\FG{\SL_n}=\{\id\}$, and $\Wx=\Waf$ is the usual affine Weyl group.
\item For the adjoint group $G=\PGL_n(\C)$, $\Is=\Iaf^s=\Iaf=\Z/n\Z$ and $\FG{\PGL_n}\cong \Z/n\Z$ acts on $\Iaf$ such that $\pi_i=\pi_{i+n\Z}$ rotates the affine Dynkin diagram forward by $i$ nodes. The group $\Wx$ is the usual extended affine Weyl group.
Since $\pi_i^{-1}=\pi_{-i}$, we have $i^*=-i\in \Z/n\Z.$
\item 
For $G=\GL_n(\C)$, there are group isomorphisms $X^\vee/Q^\vee\cong \FG{\GL_n} \cong \Z$. Let $\sh$ denote the generator of $\FG{\GL_n}$. 
The action of $\FG{\GL_n}$ on the affine Weyl group $\Waf$ is given by the
map $\FG{\GL_n}\to \Aut(\Iaf)^s$, which sends $\sh$ to $\pi_1=\pi_{1+n\Z}$.
Using the notation for $\eps_i$ in \S \ref{SS:type A root data}
Equation \eqref{E:Pi as translation finite} reads
\begin{align*}
\sh=t_{\varepsilon_1} u_1
\end{align*}
where $u_1 = s_1s_2\dotsm s_{n-1}\in\Wfin.$ 
One may verify that $\sh\, s_i\, \sh^{-1} = s_{i+1}$ with indices taken modulo $n\Z$.
Let $\eps=\sum_{i=1}^n \eps_i \in X^\vee(GL_n)$. We have
\begin{align}\label{E:Pi GL as translation finite}
  \sh^{rn+i} = t_{\eps}^r t_{\eps_1+\dotsm+\eps_i} u_1^i \qquad\text{for $r\in\Z$ and $0\le i\le n-1$.}
\end{align}
\end{enumerate}
\end{ex}

An element $\fundel\in \FG{G}$ acts on $X$ by $\cl(\fundel)$.

For $i\in I$ we denote by $P_i$ 
the standard maximal parabolic subgroup of $G$ corresponding to $i$; 
$P_i$ is the subgroup containing $B$
such that its Weyl group is 
$W_{P_i}=\langle s_j\mid j\in I\setminus \{i\}\rangle.$
Let $W_G^{P_i}$ denote the set of minimum length representatives of the coset space $\Wfin/W_{P_i}$.
\begin{ex}
For $G=\PGL_n(\C)$, we have
$\Is=\Iaf\cong \Z/n\Z$.
The element $u_i\;(1\le i\le n-1)$ is 
the unique  $i^*$-Grassmannian
element, i.e. an element of $W_G^{P_{i^*}}$,  of the maximal length $i(n-i)$
in $\Wfin=S_n$. 
For example, for $n=5$
\[
u_1=
s_1s_2s_3s_4,\quad
u_2=s_2s_3s_4s_1s_2s_3,\quad
u_3=s_3s_4s_2s_3s_1s_2,\quad
u_4=s_4s_3s_2s_1.\]
\end{ex}
\begin{ex}\label{Ex:type C un}
For the adjoint group of type $C_n$,
$\Is=\{0,n\}$. Note that $n^*=n.$ Then $u_n$
is the unique element in 
$W_G^{P_n}$ having the maximal length $n(n+1)/2$. 
Explicitly,
\begin{equation*} u_n=s_n(s_{n-1}s_n )\cdots(s_2\cdots s_{n-1}s_n)(s_1\cdots s_{n-1}s_n).
\end{equation*}
We have $\Wfin\cong S_n\ltimes \{\pm 1\}^n$ with simple generators
$s_1,\dotsc,s_{n-1}$ generating $S_n$ and $s_n$ mapping to the 
element $(\id, (1,1,\dots,1,-1))$.
There is a standard isomorphism of $\Wfin$ with signed permutations,
the subgroup of permutations $w$ of the set $\{\pm1,\pm2,\dotsc,\pm n\}$
such that $w(-i)=- w(i)$ for $1\le i\le n$. Then $s_1,\dotsc,s_{n-1}$ are defined
by their actions on $\{1,2,\dotsc,n\}$ and $s_n$ fixes $1,2,\dotsc,n-1$ and 
$s_n(n)=-n$. The element $u_n$ is defined by $u_n(j)=-(n+1-j)$ for $1\le j\le n$.
\end{ex}

\subsubsection{Bruhat order, length function, and affine Grassmann elements for $\Wx$}
Let $\le$ denote the Bruhat order on the Coxeter group $\Waf$.
The partial order $\le$ extends to $\Wx$: $\fundel v \le \fundel' w$ with $\fundel,\fundel'\in\FG{G}$
and $v,w\in\Waf$ if $\fundel=\fundel'$ and $v\le w$.
The length function extends to $\Wx$ by $\ell(\fundel v)=\ell(v)$ for all $\fundel\in \FG{G}$ and $v\in\Waf$. This agrees with counting inversions (positive affine real roots sent to negative) since elements of $\FG{G}$ permute the affine simple roots and therefore preserve the cone of affine positive real roots. 
An element $w\in \Wx$ is said to be affine Grassmannian (denoted $w\in \Wxgr$) if
$ws_i>w$ for all $i\in I$. Let $\Wafgr = \Waf\cap \Wxgr$.
We have $\Wxgr=\FG{G}\Wafgr$.

\subsection{$K$-theoretic nil-Hecke algebra}
\label{SS:K nilHecke}

We recall the construction in \cite{LSS:K} of a $K$-theoretic analogue of Peterson's affine level zero nil-Hecke algebra \cite{Lam}, \cite{Pet}, \cite{LS:Q=Aff}.
We introduce a version $\Kx$ that uses the extended affine Weyl group $\Wx$ rather than the affine Weyl group; see \cite{CMP}, \cite{LS:double Toda} for the cohomological analogue of this extended version.
Our construction in this article includes nonsemisimple groups such as $\GL_n(\C)$.

\subsubsection{Definition of $\Kx$ for reductive $G$}\label{SS: def Kx}

Let $\Rep\cong \Z[X]=\bigoplus_{\la\in X} \Z e^\la$ be the representation ring of $T$.
Since $\Wfin$ acts naturally on $X$ it acts on $\Rep$.
Define the map $c: X\to \Rep$ by
\begin{align*}
c(\al)=1-e^{\al}\qquad\text{for $\al\in X$.}
\end{align*}
Note that
\begin{align}\label{E:c negative}
c(-\alpha)=-e^{-\alpha} c(\alpha).
\end{align}
Since $\Wfin$ acts on $\Rep$ and $\Phi$,
it also acts on the localization $\Reprat$ of $\Rep$ at 
the multiplicatively closed set $\Delta$ generated by the elements $c(\al)$
for $\al\in \Phi$. Hence $\Wx$ acts on $\Rep$ and $\Reprat$ by the level zero action.

Let $\Kx^\Delta$ be the smash product of $\Reprat$ and $\Z[\Wx]$,
which is by definition the left $\Reprat$-module $\Kxrat = \Reprat\otimes_\Z \Z[\Wx]$ equipped with the product
\begin{equation}\label{eq:Kxrat product}
(f\otimes w)(g\otimes v)=fw(g)\otimes wv \quad \text{for $f,g\in \Reprat$}
\end{equation}
and $w,v\in \Wx$. $\Kxrat$ acts on $\Reprat$, $\Wx$ by the level zero action
and $\Reprat$ by left multiplication.

For $i\in \Iaf$, define $T_i,D_i\in \Kxrat$ by
\begin{align*}
T_i&=c(\al_i)^{-1}(s_i-1),\\
D_i&=T_i + 1.
\end{align*}
The elements $T_i$ satisfy $T_i^2=-T_i$ and the braid relations. 
The elements $D_i$ satisfy $D_i^2=D_i$ and the braid relations.
Therefore it makes sense to define $T_w$ and $D_w$ for $w\in \Waf$.

We have \cite[(2.6)]{LSS:K}
\begin{align}\label{E:T derivation}
  T_i f &= T_i(f) + s_i(f) T_i\qquad\text{for $i\in \Iaf$ and $f\in R(T)$.}
\end{align}

For $\fundel\in\FG{G}$ and $w\in\Waf$, we define $T_{\fundel w}, D_{\fundel w}\in \Kx$ by
\begin{align}
\label{E:T pi Waf}
  T_{\fundel w} &= \fundel T_w, \\
\label{E:D pi Waf}
  D_{\fundel w} &= \fundel D_w.
\end{align}
In particular 
\begin{align*}
    T_\fundel=\fundel=D_\fundel\qquad\text{for $\fundel\in\FG{G}$.}
\end{align*}
By \eqref{E:pi on simple reflections} and \eqref{E:pi on simple roots}, 
we deduce that for all $\fundel\in \FG{G}$ and $i\in \Iaf$ we have
\begin{align}\label{E:pi conj T}
\fundel T_i \fundel^{-1} &= T_{\fundel(i)}, \\
\label{E:pi conj D}
\fundel D_i \fundel^{-1} &= D_{\fundel(i)}.
\end{align}
From \eqref{eq:Kxrat product}, we have
\begin{align}\label{E:pi conj R(T)}
    \fundel f \fundel^{-1} = \fundel(f)\qquad\text{for $f\in R(T)$ and $\fundel\in \FG{G}$.}
\end{align}
Define the affine level zero $K$-nil-Hecke ring $\Kx$ to be the
subring of $\Kxrat$ generated by $R(T)$, $\{T_i\mid i\in \Iaf\}$, and $\FG{G}$.
By \eqref{E:T derivation} and \eqref{E:pi conj T}
it follows that $\Kx$ acts on $\Rep$ and satisfies 
\begin{align}\label{E:T basis}
\Kx=\bigoplus_{w\in \Wx} \Rep T_w.
\end{align}

\begin{lem} For $w_1,w_2\in\Wx$ such that $\ell(w_1)+\ell(w_2)=\ell(w_1w_2)$, we have
\begin{align}
\label{E:TT}
T_{w_1} T_{w_2} &= T_{w_1w_2}, \\
\label{E:DD}
D_{w_1} D_{w_2} &= D_{w_1w_2}.
\end{align}
\end{lem}
\begin{proof} Let $w_1=\fundel_1 v_1$ and $w_2=\fundel_2 v_2$ where $\fundel_1,\fundel_2\in \FG{G}$ and
$v_1,v_2\in\Waf$. Then $w_1w_2=\fundel_1 v_1 \fundel_2 v_2 = \fundel_1 \fundel_2 (\fundel_2^{-1}(v_1))v_2$.
So $T_{w_1w_2} = \fundel_1\fundel_2 T_{\fundel_2^{-1}(v_1) v_2}$ by definition.
Since $\FG{G}$ acts by Dynkin automorphisms on $\Waf$ we have
$\ell(\fundel_2^{-1}(v_1))=\ell(v_1)$. 
It follows that 
$$\ell(\fundel_2^{-1}(v_1))+\ell(v_2)=\ell(w_1)+\ell(w_2)=\ell(w_1w_2)=\ell(\fundel_2^{-1}(v_1) v_2).$$ 
Therefore $T_{\fundel_2^{-1}(v_1) v_2}=T_{\fundel_2^{-1}(v_1)} T_{v_2}$. Then 
\begin{align}
    T_{w_1w_2} 
    &= \fundel_1\fundel_2 T_{\fundel_2^{-1}(v_1)} T_{v_2} && 
    \\
    &= \fundel_1 T_{v_1} \fundel_2 T_{v_2} && \text{by \eqref{E:pi conj T}} \\
    &= T_{w_1} T_{w_2} &&\text{by \eqref{E:T pi Waf}}
\end{align}
as required. The proof of \eqref{E:DD} is the same (replace $T$ by $D$).
\end{proof}

\begin{lem} \label{L:Ti on Tw}
For all $i\in \Iaf$ and $w\in\Wx$
\begin{align}
\label{E:Ti on Tw}
    T_i T_w &= \begin{cases}
        T_{s_i w} & \text{if $s_i w > w$} \\
        - T_w &\text{if $s_i w<w$,}
    \end{cases} \\
\label{E:Di on Dw}
D_i D_w &= \begin{cases}
    D_{s_iw} & \text{if $s_iw>w$} \\
    D_w &\text{otherwise.}
\end{cases}
\end{align}
\end{lem}
\begin{proof} In the case $s_iw>w$ \eqref{E:Ti on Tw}
and \eqref{E:Di on Dw} hold by definition.
Suppose $s_iw<w$. Again by definition we have
$T_w = T_i T_{s_i w}$ and $D_w = D_i D_{s_iw}$.
We have $T_i T_w=T_i^2 T_{s_i w} = -T_i T_{s_iw}=-T_w$
and $D_i D_w = D_i^2 D_{s_i w} = D_i D_{s_i w} = D_w$.
\end{proof}

\subsubsection{$\Kx$ as an extended version of the simply-connected construction}
Consider the groups $G$ and $\Gsc$.
The ring $\mathbb{K}_{\Gsc}$ is the original $K$-nil-Hecke ring defined in \cite{LSS:K}.
The ring $\Kx$ may be viewed as an ``extended" version of $\K$ as follows.

First, there is a natural embedding of rings $\K\to \Kx$.
There is an injective ring map $R(\Tsc) \cong \Z[Q] \subset \Z[X] \cong R(T)$.
Since $c(\alpha_i)\in R(\Tsc)$, $R(\Tsc)^\Delta \subset R(T)^\Delta$
and $T_i \in \K$ for all $i\in\Iaf$. Finally $\FG{\Gsc}$ is trivial. 
This shows that there is an inclusion of rings $\K \subset \Kx$.
Moreover \eqref{E:pi conj T} and \eqref{E:pi conj R(T)} 
show that the group $\FG{G}$ acts on
$\K\subset \Kx$ by conjugation since the action of $\Wfin$ on $X$ stabilizes $Q$.
It is not hard to check that $\Kx$ is isomorphic to the smash product $\Z[\FG{G}] \otimes_\Z R(T) \otimes_{R(\Tsc)} \K$ 
with multiplication defined by 
$(\fundel\otimes a)(\fundel'\otimes b) = \fundel \fundel' \otimes (\fundel')^{-1}(a)b$ for $\fundel,\fundel'\in\FG{G}$ and $a,b\in R(T) \otimes_{R(\Tsc)} \K$. We usually omit the tensor symbol and write
\begin{align}\label{E:Pi on K}
  \fundel a \fundel^{-1} = \pi(a)\qquad\text{for $\fundel\in\FG{G}$ and $a\in R(T) \otimes_{R(\Tsc)} \K$.}
\end{align}
$\Kx$ is a $\K$-$\K$-bimodule, graded by the finite group $\FG{G}$. 
It naturally contains $\K$ as a subring.

\subsubsection{Relation between $T_w$ and $D_w$}
\begin{prop} \label{P:T and D}
For all $w\in \Wx$
\begin{align}\label{E:D into T}
D_w = \sum_{v\le w} T_v
\end{align}
\end{prop}
\begin{proof} 
We first prove \eqref{E:D into T} for $w\in \Waf$. 
We use induction on $l=\ell(w)$. For $l=0$ we have $w=\id$
and $D_\id = \id = T_\id$, proving \eqref{E:D into T} in this case.
Suppose $l\ge 1$ and \eqref{E:D into T} holds for $u\in \Waf$ for $\ell(u)<l$. Let $w\in \Waf$ be an element of length $l$. There are $i\in \Iaf$ and $u\in \Waf$, $\ell(u)=l-1$
such that $w=s_iu$. We have
\begin{align*}
  D_w &= D_i D_u= D_i \sum_{x\le u} T_x \\
  &= (1 + T_i) \sum_{x\le u} T_x \\
  &= \sum_{x\le u} T_x - \sum_{\substack{x\le u \\ s_i x < x}} T_x + \sum_{\substack{x\le u \\ s_i x > x}} T_{s_i x}\quad\text{by Lemma \ref{L:Ti on Tw}} \\
  &= \sum_{\substack{x\le u \\ s_i x > x}} (T_x + T_{s_i x}) \\
  &= \sum_{y\le w} T_y.
\end{align*}
The last equality uses the Diamond Lemma (see for example \cite[Proposition 5.4.3]{MP}) for the Bruhat order (see also \cite[Appendix A]{IIN}).

In general we write $w'\in\Wx$ as $\fundel w$ for $\fundel\in \FG{G}$ and $w\in \Waf$.
We have
\begin{align*}
    D_{w'} = \fundel D_w = \fundel \sum_{v\le w} T_v 
    = \sum_{v\le w} T_{\fundel v} = \sum_{v'\le w'} T_{v'}.
\end{align*}
\end{proof}

\subsubsection{Some elements of $\Kx$}

We view $\Wx \subset \Kx$ by rewriting the definition of $T_i$:
\begin{align}\label{E:s in K}
  s_i = 1 + c(\al_i) T_i\qquad\text{for $i \in \Iaf$.}
\end{align}
For $\alpha\in\Phi$ define $T_\alpha\in \Kxrat$ by
\begin{align*}
    T_\alpha = c(\alpha)^{-1}(s_{\alpha}-1).
\end{align*}
Then we have
\begin{equation}
\label{E: s_alpha in K}
s_\alpha=1+   c(\al) T_\alpha.  
\end{equation}
Let $w\in \Wfin$ and $i\in I$ such that $\alpha=w(\alpha_i).$ Then
\begin{equation}
T_\alpha=w T_i w^{-1}.
\end{equation}

\begin{lem} \label{L:T alpha} 
$T_\alpha\in \Kx$ for $\alpha\in \Phi$.
\end{lem}
\begin{proof}
Let $i\in I$ and $v\in \Wfin$ be such that
$v(\alpha_i)=\alpha$. Then by definition $v(\alpha_i^\vee)=\alpha^\vee$.
We have $T_\alpha = v T_i v^{-1}$. Since $\Wfin\subset \Kx$ via \eqref{E:s in K}
it follows that $T_\alpha\in \Kx$.
\end{proof}

\begin{lem}\label{L:T0 on 1}
  $c(\theta)^{-1} (t_{\theta^\vee}-1)\in \Kx$.  
\end{lem}
\begin{proof}
Using \eqref{E: s_alpha in K} and \eqref{E:c negative}, we have 
\begin{align*}
t_{\theta^\vee}-1&=s_0s_\theta-1\\
&=(1+c(-\theta)T_0)(1+c(\theta)T_\theta)-1\\
&=c(-\theta)T_0
+c(\theta)T_\theta
+c(-\theta)T_0
c(\theta)T_\theta\\
&=c(\theta)
\left(
-e^{-\theta}T_0
+T_\theta
-e^{-\theta}T_0
c(\theta)T_\theta
\right)
\end{align*}
as required.
\end{proof}

\subsection{$K$-Peterson subalgebra and the Schubert bases of $K^T_*(\Gr_G)$}
\label{SS:Peterson subalgebra}

The extended $K$-theoretic Peterson algebra is by definition the
centralizer $\Lx= Z_{\Kx}(\Rep$).
Let $\Lxrat = \Reprat \otimes_{\Rep} \Lx$.

\subsubsection{The $K$-theoretic Peterson algebra $\mathbb{L}_G$} The same proof for the simply-connected case in \cite[Lemma 5.2]{LSS:K} works
for reductive $G$.

\begin{lem}\label{L:L rat}
$\Lxrat$ has an $\Reprat$-basis given by $t_\la$ for $\la\in X^\vee$.
Equivalently we have a natural $\Rep^\reg$-algebra isomorphism 
\begin{align}
\label{E:L-reg}
\Lx^\reg &=\Rep^\reg[X^\vee].
\end{align}
Moreover
\begin{align*}
\Lx &= \Lxrat \cap \Kx.
\end{align*}
\end{lem}

$\Lx$ is a $\Rep$-Hopf algebra. 
Over $\Reprat$, the coproduct is a $\Reprat$-algebra homomorphism
defined by $\Delta(t_\la)= t_\la\otimes t_\la$ for $\la \in X^\vee$.
One may show that this restricts to a $\Rep$-algebra homomorphism $\Lx\to \Lx\otimes_{\Rep} \Lx$.
$\Lx$ has antipode $S(t_\la)=t_{-\la}$
and counit $\epsilon(t_\la)=\delta_{\la,0}$ for $\la\in X^\vee$.

The proof of \cite[Thm. 5.4]{LSS:K} for the simply-connected case, works in the same way for $G$ reductive. It can be used to show that for $w\in\Wxgr$ there is a unique element $k_w\in \Lx$ of the form
\begin{align*}
k_w=T_w+\sum_{v\in \Wx\setminus \Wxgr}
k_w^v T_v
\end{align*}
with $k_w^v\in \Rep$.  
Let $\ell_w\in \Lx$ be defined by
\begin{align}\label{E:l def}
  \ell_w = \sum_{\substack{v\in\Wxgr\\ v\le w}} k_v.
\end{align}

Let $K_T^*(\bGr_G)$ be the $T$-equivariant $K$-cohomology of the thick affine Grassmannian
$\bGr_G$ \cite{Ka},\cite{KaS}. There is an $R(T)$-bilinear pairing \cite{Ku} 
$\pair{\cdot}{\cdot}: K^T_*(\Gr) \otimes_{R(T)} K_T^*(\bGr_G) \to R(T)$.
For $\la\in X^\vee$ let $i_{t_\la}^*\in K^T_*(\Gr_G)$ be such that
for $\gamma\in K_T^*(\bGr_G)$, $\pair{i_{t_\la}^*}{\gamma}$ is the localization of $\gamma$ at the $T$-fixed point $t_{\la}$.

\begin{thm}\label{T:G iso}
\cite[Thm. 5.3]{LSS:K}
There is a $\Rep$-Hopf algebra isomorphism
\begin{align}
\label{E:G iso}
\Lx &\cong K_*^T(\Gr_G)
\end{align}
defined over $\Reprat$ by
\begin{align*}
  t_\la \mapsto i_{t_\la}^*.
\end{align*}
Moreover we have
\begin{align*}
k_w &\mapsto \cI_{w}\qquad\text{for $w\in \Wxgr$}, \\
\ell_w&\mapsto \O_{w}\qquad\text{for $w\in \Wxgr$},
\end{align*}
where $\cI_{w}$ is the ideal sheaf
of the boundary $\partial X_w$ of the Schubert variety $X_w$ in the affine
Grassmannian $\Gr_G$, and 
$\O_{w}$ is the structure sheaf of the
Schubert variety $X_w$ in $\Gr_G$.
\end{thm}

\begin{rem} \label{R:KGr}
Using localization \cite[Thm. 5.4]{LSS:K} shows that the basis $k_w$ of $\Lx$ corresponds to the basis $\xi_w$ of $K^T_*(\Gr_G)$ that is dual (under an evaluation pairing) 
to the structure sheaf basis of the equivariant $K$-cohomology $K_T^*(\bGr_G)$ 
of Kashiwara's thick affine Grassmannian $\bGr_G$ \cite{Ka} \cite{KaS}.
In \cite{LLMS} the above duality is interpreted via Kumar's tensor pairing \cite[\S 3]{Ku} $K_T^*(\bGr_G) \otimes_{\Rep} K^T_*(\Gr_G)\to\Rep$, yielding the equality $\xi_w=\cI_w$ in $K^T_*(\Gr_G)$. Then it is deduced in \cite[Thm. 1]{LLMS} that $\ell_w$ and $\O_w$ correspond.
\end{rem}

\subsubsection{The map $\kappa$}

Let $\kappa: \Kxrat\to\Lxrat$ be the left $\Reprat$-module homomorphism defined by
\begin{align*}
\kappa(t_\la w) &= t_\la\qquad\text{for $\la\in X^\vee$ and $w\in \Wfin$.}
\end{align*}

\begin{lem} \label{L:kappa} For all $a,a'\in \Kxrat$, $b\in \Lxrat$, $u\in \Wfin$,
and $i\in I$
\begin{align}
\label{E:kappa fixes L}
    \kappa(b) &= b, \\
\label{E:kappa Wfin}
\kappa(au) &= \kappa(a), \\
\label{E:kappa kills}
\kappa(aT_i) &= 0,\\
\label{E:kappa K action}
    \kappa(a'a) &= \kappa(a'\kappa(a)), \\
\label{E:kappa L linear}
    \kappa(b\,a) &= b\, \kappa(a).
\end{align}
\end{lem}
\begin{proof} 
Equation \eqref{E:kappa fixes L} holds since $\Lxrat$ has $\Reprat$-basis given by translation elements. For \eqref{E:kappa Wfin} by linearity let $a=t_\la w$ for $\la\in X^\vee$ and $w\in\Wfin$
so that $\kappa(au)=\kappa(t_\la wu)=t_\la=\kappa(t_\la w)=\kappa(a)$.
For \eqref{E:kappa kills}, using \eqref{E:kappa Wfin} at the end we have
\begin{align*}
\kappa(a T_i) &= \kappa(a c(\alpha_i)^{-1} (s_i - 1)) \\
&=\kappa(a c(\alpha_i)^{-1} s_i) - \kappa(a c(\alpha_i)^{-1}) = 0.
\end{align*}
For \eqref{E:kappa K action} let $a'=t_\la u$ and $a=t_\mu v$
for $\la,\mu\in X^\vee$ and $u,v\in \Wfin$.
Since $a'a=t_{\la + u(\mu)} uv=a'\kappa(a)v$, we obtain \eqref{E:kappa K action} by \eqref{E:kappa Wfin}.

For \eqref{E:kappa L linear} we have
$\kappa(ba)=\kappa(b \kappa(a))=b \kappa(a)$ by \eqref{E:kappa K action}
and \eqref{E:kappa fixes L}.
For \eqref{E:kappa L linear}, we may assume $b=t_\la$ for $\la\in X^\vee$. By linearity let $a=t_\mu u$ for $\mu\in X^\vee$ and $u\in \Wfin.$ We have
$\kappa(ba)=\kappa(t_{\la+\mu}u)=t_{\la+\mu}=t_\la \kappa(a)=b\kappa(a)$ as required.
\end{proof}

\begin{lem} \label{L:kappa kernel T}
For $w\in \Wx\setminus \Wx^0$ and $a\in \Kx$
\[
\kappa(aT_w)=0.
\]
\end{lem}
\begin{proof}
    This follows from \eqref{E:kappa kills}.
\end{proof}

The following was proved in \cite{LLMS} for $\K$ but the same proof
works for $\Kx$.

\begin{prop} \label{P:kappa basis}
\cite[Prop. 2]{LLMS}
The map $\kappa$ restricts to a left $\Rep$-module homomorphism
$\Kx\to\Lx$ such that
\begin{align}
\label{E:kappa basis}
\kappa(T_w) &= k_w &&\text{for $w\in \Wx^0$}, \\
\label{E:kappa l}
\kappa(D_w) &= \ell_w &&\text{for $w\in \Wx^0$.}
\end{align}

\end{prop}

\begin{rem} $\kappa$ also restricts to a left $\Rep$-module homomorphism
$\K\to\LL$.
\end{rem}

\subsection{The ``star" action of $\Kx$ on $\Lx$}
\label{SS:star action}

\begin{lem}\label{L:star action}
There is an action $*$ of $\Kxrat$ on $\Lxrat$ extending the 
left $\Reprat$-module structure, defined on an element $b\in\Lxrat$ by
\begin{align*}
   w * b &= w b w^{-1} &\qquad&\text{for $w\in \Wfin$} \\
   t_\la* b &= t_\la b && \text{for $\la\in X^\vee$ and $b\in \Lxrat$.}
\end{align*}
Moreover this restricts to an action of $\Kx$ on $\Lx$: for $b\in \Lx$ we have
\begin{align*}
  T_i * b &= T_i b + b\,T_i + T_i \,b\, c(\al_i) T_i &&\text{for $i\in I$,} \\
  T_0 * b &= t_{\theta^\vee} c(\al_0)^{-1} (s_\theta-1)* b + c(\al_0)^{-1}(t_{\theta^\vee}-1) b. 
\end{align*}
\end{lem}
\begin{proof} We first show that 
\begin{align}\label{E:Wx stabilizes Lx}
\Wx * \Lx \subset\Lx.
\end{align}
Let $\la\in X^\vee$, $w\in \Wfin$, and $b\in\Lx$. 
For all $r\in\Rep$ we have 
\begin{align*}
  r(t_\la w * b) &= t_\la r w b w^{-1} \\
  &= t_\la w (w^{-1}(r)) b w^{-1} &\quad& \text{by 
  \eqref{eq:Kxrat product}}\\
  &=t_\la  w b (w^{-1}(r)) w^{-1} &&\text{since $b\in \Lx$}\\
  &= t_\la w b w^{-1} r &&\text{by \eqref{eq:Kxrat product}}\\
  &= ((t_\la w) *b) r.
\end{align*}
This proves \eqref{E:Wx stabilizes Lx}. We deduce that
$\Kxrat*\Lxrat\subset\Lxrat$ since $\Wx$ is a left $\Reprat$-basis of $\Kxrat$.

Next we verify the action axiom
\begin{align}\label{E:star action property}
  (aa')*b &= a*(a'*b) \qquad\text{for $a,a'\in\Kxrat$ and $b\in\Lxrat$.}
\end{align}
Let $\la,\mu\in X^\vee$ and $w,v\in \Wfin$ and $b\in\Lxrat$.
We have
\begin{align*}
(t_\la w)* ((t_\mu v) * b) &= (t_\la w) ( t_\mu v b v^{-1}) \\
&= t_\la w  t_\mu v b v^{-1} w^{-1} \\
&= t_{\la + w(\mu)} wv b (wv)^{-1} \\
&= (t_{\la+w(\mu)} wv)* b \\
&= ((t_\la w)(t_\mu v))*b .
\end{align*}
Now we check that $\Kx * \Lx \subset \Lx$. The ring $\Kx$ is generated by $\Rep$, $\{T_i\mid i\in \Iaf\}$, 
and $\FG{G}$. Since $\Lx$ is a $\Rep$-algebra and by \eqref{E:Wx stabilizes Lx},
it is enough to show that $T_i*\Lx\subset\Lx$ for $i\in\Iaf$. By Lemma \ref{L:L rat}
it is enough to show that $T_i*\Lx\subset \Kx$. Let $b\in\Lx$ and $i\in I$. We have
\begin{align}
  T_i*b &= c(\al_i)^{-1} (s_i-1)*b \nonumber\\
  &= c(\al_i)^{-1} (s_i b s_i - b)\nonumber \\
  &= c(\al_i)^{-1} ( ( 1 + c(\al_i) T_i) b (1+c(\al_i)T_i) - b)\nonumber \\
  &= c(\al_i)^{-1}( c(\al_i)T_i \,b + b\, c(\al_i) T_i + c(\al_i) T_i\, b\, c(\al_i)T_i) \nonumber\\
  &= \label{E:T_i*a}T_i\, b + b\, T_i + T_i\, b\, c(\al_i) T_i\in \Kx.
\end{align}
For $i=0$ by Lemma \ref{L:T0 on 1} we have
\begin{align}\label{E:T0 star 1}
  T_0 * 1 = c(\al_0)^{-1}(t_{\theta^\vee} - 1)  \in \Kx.
\end{align}
For $b\in\Lx$ we compute
\begin{align*}
  T_0 * b &= c(\al_0)^{-1}(t_{\theta^\vee}s_\theta - 1) * b \nonumber\\
  &= c(\al_0)^{-1}( t_{\theta^\vee} s_\theta b s_\theta - b) \nonumber\\
  &= c(\al_0)^{-1}( t_{\theta^\vee} s_\theta b s_\theta - t_{\theta^\vee} b+ t_{\theta^\vee} b- b)\nonumber \\
  &= t_{\theta^\vee} c(\al_0)^{-1} (s_\theta-1)*b + c(\al_0)^{-1}(t_{\theta^\vee}-1)\,b \\
  &= - t_{\theta^\vee} e^\theta T_\theta * b + (T_0*1) b.
  \nonumber 
\end{align*}
To show $T_0*b \in\Kx$ it suffices to show that $T_\theta*b\in \Kx$ and
$T_0*1\in\Kx$. The latter holds by \eqref{E:T0 star 1}.
For the former, there is a $u\in\Wfin$ and $i\in I$ such that $T_\theta = u T_i u^{-1}$. 
Then $T_\theta * b = u T_i u^{-1} * b \in \Kx$ by 
\eqref{E:T_i*a} and \eqref{E:Wx stabilizes Lx}.
\end{proof}

\begin{lem}\label{L:kappa action} We have
\begin{align*}
    a * b = \kappa(ab) \qquad\text{for $a\in\Kx$ and $b\in\Lx$.}
\end{align*}
\end{lem}
\begin{proof} By linearity and working over $\Reprat$ we may assume
$a = t_\mu v$ and $b=t_\la$ for $\mu,\la\in X^\vee$ and $v\in \Wfin$. We have
\begin{align*}
    \kappa(ab) &= \kappa(t_\mu v t_\la) \\
    &= \kappa(t_\mu t_{v(\la)} v)  \\
    &= \kappa(t_\mu t_{v(\la)}) &\qquad&\text{by \eqref{E:kappa Wfin}}\\
    &= t_\mu t_{v(\la)} && \text{by \eqref{E:kappa fixes L}} \\
    &= t_\mu v * t_\la \\
    &= a*b.
\end{align*}
\end{proof}

\begin{prop} \label{P:star basis formula}
For any $w\in \Wxgr$
\begin{align}
\label{E:star k basis}
    k_w &= T_w * 1, \\
\label{E:star l basis}
    \ell_w &= D_w * 1.
\end{align}
In particular $\Lx= \Kx*1$.
\end{prop}
\begin{proof}
We have $k_w = \kappa(T_w) = T_w * 1$
by \eqref{E:kappa basis} and Lemma \ref{L:kappa action}.
This proves \eqref{E:star k basis}.
Note that if $v\in\Wx\setminus \Wxgr$ then $T_v * 1 = 0$ since $T_i*1=0$ for $i\in I$.
For \eqref{E:star l basis} we have
\begin{align*}
  D_w*1 = \sum_{v\le w} T_v * 1 = \sum_{v\le w, \, v\in \Wxgr} k_v = \ell_w
\end{align*}
by Prop. \ref{P:T and D}, \eqref{E:star k basis}, and \eqref{E:l def}.
\end{proof}

\begin{prop}\label{P:T on k}
For $i\in \Iaf$ and $w\in \Wx^0$
\begin{align}
\label{E:T on k}
  T_i * k_w &= 
  \begin{cases}
      k_{s_i w} & \text{if $s_iw>w$ and $s_iw\in \Wx^0$} \\
      0&\text{if $s_iw>w$ and $s_iw\not\in \Wx^0$} \\
      -k_w &\text{if $s_iw<w$}
  \end{cases},\\
\label{E:D on l}
  D_i * \ell_w &= 
  \begin{cases}
      \ell_{s_i w} & \text{if $s_iw>w$ and $s_iw\in \Wx^0$} \\
      \ell_w &\text{otherwise}
  \end{cases}.
\end{align}
\end{prop}
\begin{proof} We have
\begin{align*}
  T_i * k_w &= T_i * (T_w * 1) &\qquad&\text{by \eqref{E:star k basis}} \\
  &= (T_i T_w) * 1 &&\text{by \eqref{E:star action property}} \\
  &= \kappa(T_i T_w) &&\text{by Lemma \ref{L:kappa action}.} 
\end{align*}
If $s_iw <w$ then $s_iw\in\Wxgr$ and $\kappa(T_iT_w)=\kappa(-T_w)=-k_w$
by \eqref{E:Ti on Tw} and \eqref{E:kappa basis}.
Otherwise $s_iw>w$ and $\kappa(T_i T_w)=\kappa(T_{s_iw})$ by \eqref{E:Ti on Tw}.
If $s_iw\notin\Wxgr$ then $\kappa(T_{s_iw})=0$ by Lemma \ref{L:kappa kernel T}.
If $s_iw\in\Wxgr$ then $\kappa(T_{s_iw})=k_{s_iw}$ by \eqref{E:kappa basis}.
This proves \eqref{E:T on k}.

The proof of \eqref{E:D on l} is similar.
\end{proof}

\begin{cor}\label{C:TD on kl}
For any $v,w\in\Wx$ such that
$vw\in\Wxgr$ and $\ell(v)+\ell(w)=\ell(vw)$,
we have $w\in\Wxgr$ and
\begin{align}
\label{E:Ts on k}
k_{vw} &= T_v * k_w \\
\label{E:Ds on l}
\ell_{vw} &= D_v * \ell_w.
\end{align}
\end{cor}
\begin{proof} The fact that $w\in\Wxgr$ is standard.
We have $k_{vw} = T_{vw}*1= T_v*(T_w*1)=T_v*k_w$
using \eqref{E:star k basis}, \eqref{E:TT}, \eqref{E:star k basis}.
The same proof works for \eqref{E:Ds on l}.
\end{proof}

Recall the notation $\fundel(a)=\fundel a\fundel^{-1}$ for $\fundel\in \FG{G}$ and $a\in \Kx$. 
For $\fundel=t_\fundel u_\fundel \in \FG{G}$ and $a\in \Lx$,  we note that 
\begin{equation}\label{E: pi acts on L by conjugation of finite part}
\fundel(a)=u_\fundel *a=u_\fundel au_\fundel^{-1}
\end{equation}
because $u_\fundel* a\in \Lx$ commutes with $t_\fundel$.
\begin{lem}\label{E:Pi on L} Let $\fundel\in\FG{G}$. Then
\begin{align*}
  \fundel * a = t_\fundel \fundel(a)\qquad\text{for $a\in \Lx$.}
\end{align*}
In particular for $w\in\Wxgr$
\begin{align}
\label{E:pi star on k}
    \fundel * k_w &= t_\fundel \fundel(k_w) = k_{\fundel w}, \\
\label{E:pi star on l}    
    \fundel * \ell_w &= t_\fundel \fundel(\ell_w) = \ell_{\fundel w}.
\end{align}
\end{lem}
\begin{proof} 
By the definition of the $*$ action and \eqref{E: pi acts on L by conjugation of finite part}, we have
 \begin{equation*}
   \fundel*a = t_\fundel u_\fundel a u_\fundel^{-1} 
   = t_\fundel \fundel(a).
 \end{equation*}
This gives the first equality in \eqref{E:pi star on k}. We have
\begin{align*}
\fundel * k_w &= T_\fundel * k_w = k_{\fundel w}
\end{align*}
by \eqref{E:T on k}. This proves the second equality in \eqref{E:pi star on k}.
Since $v\le w$ if and only if $\fundel v\le \fundel w$
for all $v,w\in\Wx$, \eqref{E:pi star on l} follows from \eqref{E:pi star on k}.
\end{proof}

\begin{cor} \label{C:k fundamental}
For all $\fundel\in\FG{G}$, 
\begin{align}\label{E:fundamental basis}
    k_{\fundel} = \ell_\fundel = t_{\fundel}.
\end{align}
For $G$ semisimple
\begin{align}\label{E:fundamental semisimple basis}
k_{\pi_i} = \ell_{\pi_i}=t_{\om_i^\vee}\qquad\text{for $i\in \Is$.}
\end{align}
\end{cor}
\begin{proof} We apply \eqref{E:pi star on k} for $w=\id$. By \eqref{E:fundamental basis}
$k_{\fundel} = \fundel * k_\id = \fundel * 1 = t_\fundel$. Since $\fundel$ is a Bruhat-minimal element in $\Wxgr$,
$\ell_\fundel=k_\fundel$.
\end{proof}

\subsection{Factorization and $k$-rectangle classes}
\label{SS:k rectangles}
The appropriate Schubert bases of $H^T_*(\Gr_G)$ and $K^T_*(\Gr_G)$ 
have a factorization property. The factorization phenomenon was first observed by Lapointe and Morse \cite{LM}
for $H_*(\Gr_{\SL_n})$. Magyar \cite{Mag} obtained it for equivariant homology in
general type for $G$ simply-connected.
For $K_*(\Gr_{\SL_n})$ it was proved by Takigiku \cite{Tak}.
We observe that the factorization works for reductive $G$.

\begin{lem} \label{L:finite Weyl invariants star}
Let $\Lx^{\Wfin}$ denote the $\Wfin$ invariants in $\Lx$ under the $*$ action.
Then for all $a\in \Kx$, $b\in \Lx^{\Wfin}$, and $b'\in \Lx$ we have
\begin{align}\label{E:finite Weyl invariants under star}
  a *( b\,b') = b\, (a*b')\qquad\text{for $a\in\Kx$, $b\in \Lx^{\Wfin}$ and $b'\in\Lx$.}
\end{align}
\end{lem}
\begin{proof} Working over $\Reprat$, it is enough to check \eqref{E:finite Weyl invariants under star} for $a=t_\la u$ for $\la\in X^\vee$ and $u\in\Wfin$. We have
\begin{align*} 
(t_\la u) *(b\,b') &= t_\la ubu^{-1} ub'u^{-1} 
=t_\la b ub'u^{-1}
=  t_\la b\, (u*b')= 
 b\, t_\la(u*b')
= b\, (t_\la u * b').
\end{align*}
\end{proof}

The following Lemma appears in \cite[Lemmas 3.2, 3.3]{LS:Q=Aff} for affine Weyl groups $\Waf$ but the same proofs work for extended affine Weyl groups $\Wx$ for $G$ reductive.

\begin{lem}\label{L:Grass} Let $w=ut_\la\in\Wx$ with $u\in\Wfin$ and $\la\in X^\vee$.
Then $w\in \Wxgr$ if and only if $\la$ is antidominant and for $i\in I$ such that $s_i \la=\la$,
we have $us_i > u$. Moreover, in this case we have
\begin{align*}
\ell(ut_\la) = - \langle \la, 2 \rho \rangle - \ell(u)\qquad\text{for $ut_\la\in \Wxgr$.}
\end{align*}
In particular, $t_\lambda\in \Wxgr$ for antidominant $\la\in X^\vee$.
\end{lem}

\begin{lem}\label{L:invariant} For any antidominant $\la\in X^\vee$,
$\ell_{t_\la}$ is $\Wfin$-invariant under the $*$ action.
\end{lem}
\begin{proof} 
Let $i\in I.$ If $\langle \la,\alpha_i\rangle < 0$, then 
$t_\lambda(\alpha_i)$ is a negative root 
in the affine root system. It follow that $s_i t_\lambda<t_\lambda$
(see for example \cite[(2.3.3)]{Mac_Hecke}).
If $\langle \la,\alpha_i\rangle = 0$ we have $s_i t_\la\notin \Wxgr$ by Lemma \ref{L:Grass}.
By Proposition \ref{P:T on k} we have $D_i * \ell_{t_\la} = \ell_{t_\la}$ for all $i\in I$.
But this is equivalent to $s_i *\ell_{t_\la} = \ell_{t_\la}$ for all $i\in I$.
\end{proof}

\begin{rem}  
Under the $K$-theoretic Peterson isomorphism, 
for any antidominant
element $\la\in X^\vee$, 
the Schubert basis element $\ell_{t_\la}$ corresponds to the quantum parameter 
$q_{-\la}.$
\end{rem}

\begin{prop} \label{P:k rectangle} 
Let $w= u t_\gamma\in \Wxgr$ with $u\in \Wfin$ and $\gamma\in X^\vee$, and let $\la \in X^\vee$ be antidominant.
Then $wt_\la\in\Wxgr$, $\ell(wt_\la)=\ell(w)+\ell(t_\la)$ and
\begin{align*}\label{E:k rectangle}
\ell_{w t_\la} = \ell_w \ell_{t_\la}.
\end{align*}

\end{prop}
\begin{proof} 
Let $w=ut_\la$ with $u\in\Wfin$ and $\la\in X^\vee$. Then $\la$ is antidominant
and for $i\in I$ such that $s_i(\la)=\la$ we have $us_i>u$. 
$\la+\gamma$ is antidominant. The condition $s_i(\la+\gamma)=\la+\gamma$ 
is equivalent to $0=\pair{\la+\gamma}{\alpha_i} = \pair{\la}{\alpha_i}+\pair{\gamma}{\alpha_i}$.
But both pairing values are nonpositive so they must both vanish. In particular 
it implies that $s_i(\la)=\la$. For such $i$ we have $us_i>u$. By Lemma \ref{L:Grass} we have $w t_\gamma\in\Wxgr$.
By Lemma \ref{L:Grass} we have
\begin{align*}
    \ell(w t_\gamma) &= -\pair{\la+\gamma}{2\rho} - \ell(u) \\
    &= -\pair{\la}{2\rho}-\ell(u)-\pair{\gamma}{2\rho} \\
    &= \ell(w) + \ell(t_\gamma).
\end{align*}

As for \eqref{E:k rectangle}, we have
\begin{align}
  \ell_{wt_\la} &= D_w * \ell_{t_\la} &&\text{by \eqref{E:D on l}} \\
  &= \ell_{t_\la} D_w * 1 &&\text{by \eqref{E:finite Weyl invariants under star}} \\
  &= \ell_{t_\la} \ell_w &&\text{by \eqref{E:star l basis}.}
\end{align}
\end{proof}
 Suppose $G$ is semisimple.
For $u\in \Wfin$, define the antidominant coweight
\begin{align*}
  \gamma_u = - \sum_{\substack{i\in I \\ u s_i < u}} \om_i^\vee.
\end{align*}

\begin{lem}\label{L:irr} Let $u\in \Wfin$ and $\gamma\in X^\vee$.
Then $u t_\gamma\in\Wxgr$ if and only if $\gamma-\gamma_u$ is antidominant.
\end{lem}
\begin{proof} By Lemma \ref{L:Grass} $ut_\gamma\in\Wxgr$ if and only if, for all $i\in I$,
$\pair{\gamma}{\alpha_i}\le 0$ with $\pair{\gamma}{\alpha_i}<0$ if $u s_i < u$.
Writing $\gamma= \sum_{i\in I} c_i \om_i^\vee$ we have $c_i\le0$ and
$c_i<0$ if $us_i<u$, that is, $\gamma-\gamma_u$ is antidominant.
\end{proof}

We call an element of $\Wxgr$ \textit{irreducible} if it has the form $u t_{\gamma_u}$ for $u\in \Wfin$. 

Lemma \ref{L:irr} and Proposition \ref{P:k rectangle} says that every element of $\Wxgr$ is canonically the product of a $\Wxgr$ irreducible element and an antidominant translation. 
In particular
\begin{align}\label{E:canonical k rectangle}
\ell_w = \ell_{t_{\gamma-\gamma_u}} \ell_{u t_{\gamma_u}}.
\end{align}
We shall consider another factorization that expresses the formulas using elements of $\Waf$ as much as possible.  We write translation elements in the form \eqref{E:extends affine weyl}:
\begin{align}\label{E:translation to fund aff weyl}
  t_\la &= \fundel_\la \uu_\la \qquad\text{for $\la\in X^\vee$, $\fundel_\la\in \FG{G}$, and $\uu_\la\in \Waf$.}
\end{align}
Note that the map $\la\mapsto \fundel_\la$ is a group homomorphism and that
\begin{align*}
  \uu_\la \in \Wafgr\qquad\text{if $\la$ is antidominant,}
\end{align*}
by Lemma \ref{L:irr}.

The following works for $G$ reductive but for ease of notation we assume here that
$G$ is semisimple (see Remark \ref{R: kappa m for GL} for $G=\GL_n(\C)$ case).
We define
\begin{align*}
\kappa_i = \uu_{-\om_i^\vee}\qquad\text{for $i\in I$.}
\end{align*}

\begin{prop}
We have
\begin{align}\label{E:kappa i}
    \kappa_i &= 
    \begin{cases}
        \pi_i  t_{-\om_i^\vee}=\pi_i u_i \pi_i^{-1} & \text{if $i\in\Is\setminus\{0\}$}\\
          t_{-\om_i^\vee}
          & \text{if $i\in I \setminus \Is$.}
    \end{cases}
\end{align}
\end{prop}
\begin{proof}
    For $i\in \Is\setminus\{0\}$, by \eqref{E:fundamental into translation finite} 
    $t_{-\om_i^\vee} = u_i \pi_i^{-1} = \pi_i^{-1}(\pi_i u_i \pi_i^{-1})$. We have written
    $t_{-\om_i^\vee}$ as the element $\pi_i^{-1}\in \FG{G}$ times an element
    $\pi_i u_i \pi_i^{-1}$ of $\Waf$. By the definition \eqref{E:translation to fund aff weyl}
    of $\uu_{-\om_i^\vee}$ the proposition holds. For $i\in I\setminus \Is$, $\om_i^\vee\in Q^\vee$ so that $t_{-\om_i^\vee}\in \Waf$ and $\kappa_i = t_{-\om_i^\vee}$ holds by 
    the definition \eqref{E:translation to fund aff weyl}.
\end{proof}

\begin{rem} \label{R:kappa reduced word}
Note that if $i\in \Is\setminus\{0\}$, a reduced word for $\kappa_i$ is obtained from one of $u_i$ by the group homomorphism which replaces each $s_j$ by $s_{\pi_i(j)}$.
\end{rem}

\begin{lem} \label{L:fund translation and kappa}
\begin{align}\label{E:kappa pi i}
  \ell_{t_{-\om_i^\vee}} =  
 \begin{cases} 
  t_{-\om_i^\vee} \ell_{\kappa_i} & \text{if $i\in \Is\setminus\{0\}$} \\
  \ell_{\kappa_i} & \text{if $i\in I\setminus \Is$.}
\end{cases}
\end{align}
\end{lem}
\begin{proof} There is nothing to do unless $i\in \Is\setminus\{0\}$ which we assume. 
We have
\begin{align*}
  \ell_{t_{-\om_i^\vee}} &= \ell_{\pi_i^{-1} \kappa_i} &\qquad&\text{by \eqref{E:kappa i}} \\
  &= t_{\om_{i^*}^\vee} \pi_i^{-1}(\ell_{\kappa_i}) &&\text{by \eqref{E:pi star on l} and \eqref{E:fundamental inverse into translation finite}} \\
&= t_{\om_{i^*}^\vee} (u_i^{-1}*\ell_{\kappa_i}) &&\text{by \eqref{E: pi acts on L by conjugation of finite part}} \\
  &= (u_i*t_{\om_{i^*}^\vee} )\ell_{\kappa_i} && \\
  &= t_{-\om_i^\vee} \ell_{\kappa_i}.
\end{align*}
For the fourth equality we have used Lemma \ref{L:invariant} which says that 
$\ell_{t_{-\om_i^\vee}}$ is invariant under the $*$-action by any element of $\Wfin$.
We also used that the $*$ action of $u_i\in \Wfin$ is a ring homomorphism (given by conjugation).
For the last equality, we used $u_i \om_{i^*}^\vee = - \om_i^\vee$ (Remark \ref{R:inverse in fundamental group}).
\end{proof}

\begin{lem} \label{L:rotate kappa basis}
Let $i,j\in \Is$ with $i\ne0$. We have
\begin{align*}\label{E:rotate kappa basis}
    \pi_j(\ell_{\kappa_i}) &= t_{u_j(\om_i^\vee)-\om_i^\vee} \ell_{\kappa_i}.
\end{align*}
\end{lem}
\begin{proof} We have
\begin{align*}
  \pi_j(\ell_{\kappa_i}) &= \pi_j(t_{\om_i^\vee} \ell_{t_{-\om_i^\vee}} ) &\qquad&\text{by \eqref{E:kappa pi i}} \\
  &= \pi_j(t_{\om_i^\vee}) \ell_{t_{-\om_i^\vee}}&&\text{by Lemma \ref{L:invariant}} \\
  &= t_{u_j(\om_i^\vee)} t_{-\om_i^\vee} \ell_{\kappa_i} &&\text{by \eqref{E:kappa pi i}.} 
\end{align*}
\end{proof}

\begin{rem} The translation lattice element $u_j(\om_i^\vee)-\om_i^\vee$ is
in $Q^\vee$.
\end{rem}
\begin{prop} \label{P:factor kappa}
Suppose $w\in\Wx$ and $i\in \Is\setminus\{0\}$ are such that 
$w\kappa_i\in \Wxgr$ and $\ell(w\kappa_i)=\ell(w)+\ell(\kappa_i)$.
Then $\pi_i^{-1}(w)\in \Wxgr$ and 
\begin{align*}
  \ell_{w\kappa_i} &= \ell_{\kappa_i} \pi_i(\ell_{\pi_i^{-1}(w)}).
\end{align*}
\end{prop}
\begin{proof} We have
\begin{align*}
    \ell_{w\kappa_i} &= D_w * \ell_{\kappa_i} &&\text{by \eqref{E:D on l}}\\
    &= D_w * (t_{\om_i^\vee} \ell_{t_{-\om_i^\vee}}) &&\text{by \eqref{E:kappa pi i}} \\
    &= \ell_{t_{-\om_i^\vee}} D_w * \ell_{\pi_i} &&\text{by Lemmas \ref{L:invariant} and \ref{L:finite Weyl invariants star} and \eqref{E:fundamental semisimple basis}} \\
    &= \ell_{t_{-\om_i^\vee}} \ell_{w\pi_i} &&\text{by \eqref{E:D on l}} \\
    &= \ell_{t_{-\om_i^\vee}} \ell_{\pi_i (\pi_i^{-1}(w))}   \\
    &= \ell_{t_{-\om_i^\vee}} t_{\om_i^\vee} \pi_i(\ell_{\pi_i^{-1}(w)}) &&\text{by \eqref{E:pi star on l}} \\
    &= \ell_{\kappa_i} \pi_i(\ell_{\pi_i^{-1}(w)})&&\text{by \eqref{E:kappa pi i}.}
\end{align*}
\end{proof}

\begin{ex} For type $A_5$, we have
$$u_2
\quad\text{and}\quad \kappa_2 
=(s_4s_3)(s_5s_4)(s_0s_5)(s_1s_0).
$$
The indices of simple reflections in a reduced word for $\kappa_2$ can be read from the rows of the English/matrix-style tableau of $4\times 2$ rectangular shape with box $(r,c)$ in row $r$ and column $c$ labeled by the reflection $s_{c-r}$ with indices taken modulo $6\Z$. 

\[
\begin{ytableau} s_0&s_1 \\ s_5&s_0\\  s_4&s_5\\  s_3&s_4
\end{ytableau}
\]
\end{ex}

\begin{ex}
    In type $C_3$, 
    \[
    \kappa_1=s_1s_2s_3s_2s_1s_0,\quad
    \kappa_2=(s_2s_3 s_2s_1s_0)^2,\quad
    \kappa_3=s_0s_1s_0s_2s_1s_0.
    \]
    Note that $u_3=s_3s_2s_3s_1s_2s_3$ (cf. Example \ref{Ex:type C un}).
\end{ex}

\begin{rem}\label{R: kappa m for GL}
 The construction of $\kappa_i$ for $\PGL_n(\C)$
is compatible with the structure of $\tilde{W}_{\GL_n}$. 
Let $m\in\Z$. Writing $m=rn+i$ for $r\in\Z$ and $0\le i\le n-1$, by 
\eqref{E:Pi GL as translation finite} we have $t_{\sh^m} = t_\eps^r t_{\eps_1+\dotsm+\eps_i}$ and $u_{\sh^m} = u_1^{i}$.
It follows that for $G=\PGL_n(\C)$,  $t_{-\om^\vee_i}$  is 
the image of 
$t_{\sh^m}^{-1}$ under the natural projection $\pi: \tilde{W}_{\GL_n}\to \tilde{W}_{\PGL_n}$.
Because $t_{\sh^m}^{-1}$ is an anti-dominant translation element,
we can write $t_{\sh^m}^{-1}=\sh^{-m}\kappa_m$ for some $\kappa_m\in \hat{W}_{\GL_n}^0$.
 Then 
 we have
 \begin{align*}
    \kappa_m&=\sh^{m}t_{\sh^m}^{-1}\\
    &= t_{\eps_1+\dotsm+\eps_i} u_1^{i} t_{-(\eps_1+\dotsm+\eps_i)} \\
    &= u_1^{i} t_{\eps_{n+1-i}+\dotsm+\eps_n-(\eps_1+\dotsm+\eps_i)}.
 \end{align*}
 Note that this element only depends on $i\in \Z/n \Z$. 
 Moreover, the element is equal to $\kappa_i$ as an element of $\hat{W}_{\PGL_n}^0=\hat{W}_{\GL_n}^0\subset  \tilde{W}_{\GL_n}.$
In fact, since 
the map $\pi: \tilde{W}_{\GL_n}\to \tilde{W}_{\PGL_n}$ 
 when we restrict it to $\hat{W}_{\GL_n}^0$
 is identity,
 we deduce that
 \begin{equation}
\kappa_m=\pi(\kappa_m)=\pi(\sh^m t_{\sh^m}^{-1})=
\pi(\sh^m )\pi(t_{\sh^m}^{-1})=
\pi_i t_{-\om_i^\vee}=\kappa_i,
 \end{equation}
 where used \eqref{E:kappa i} at the last equality.
 \end{rem}

For $G=\PGL_n$, every $i \in \Iaf=\Z/n\Z$ is special.
For a $k$-bounded partition $\lambda$, we denote the corresponding element in $\WSLgr$ by $x_\la$.
For $i \in \Z/n\Z$, we have
\[
\kappa_i=x_{R_i}
\] 
where $R_i=(i)^{n-i}$ (rectangle of $n-i$ rows and $i$ columns).
Suppose we have $w\in \tilde{W}_{\PGL_n}$ such that  $w\kappa_i \in \tilde{W}_{\PGL_n}^0$ and $\ell(w\kappa_i
)=\ell(w)+\ell(\kappa_i)$. 
If we write $w\kappa_i=x_\lambda$, then 
$\lambda_1=\cdots=\lambda_{n-i}=i.$
Put $\mu=(\lambda_{n-i+1},\ldots,\lambda_l)$. 
Then we write $\lambda=\mu\cup R_i$. 
We have
\[
\ell_{w\kappa_i}=\pi_i(\ell_{x_\mu})\ell_{\kappa_i}.
\]
Note that 
the $n$-residue of the upper left corner box of $\mu$ 
considered as a subset of $\lambda$ is always $i$.

\begin{ex} Let $n=3$ and 
$w=x_{2^3,1}=s_2s_0s_1s_0s_2s_1s_0\in \tilde{W}_{PGL_3}^0.$ Recall that 
\[\kappa_1=x_{R_1}=s_2s_0=\ytableaushort{0,2}*{1,1},\quad
\kappa_2=x_{R_2}=s_1s_0=\ytableaushort{01}*{2}
\]
Using Proposition 
\ref{P:k rectangle}, we can factor out elements $\ell_{t_{-\om_i^\vee}}$ . We have $w=ut_\gamma$ with $u=s_1,\gamma=-2\om_1^\vee-2\om_2^\vee$, and  $\gamma_{u}=-\om_1^\vee.$
So we have
$$
\ell_w=\ell_{-\om_1^\vee-2\om_2^\vee}\ell_{s_1t_{-\om_1^\vee}}=\ell_{-\om_1^\vee-2\om_2^\vee}\ell_{\pi_2 s_0}
=\ell_{t_{-\om_1^\vee}}\ell_{t_{-\om_2^\vee}}^2  \ell_{\pi_2 s_0}.
$$
This decomposition is depicted by 
  \begin{center}
   \setlength\unitlength{0.6truecm}
  \begin{picture}(4,5.3)(0,0)
 \put(0,0){\line(1,0){1}}
 \put(0,1){\line(1,0){1}}
 \put(0,3){\line(1,0){2}}
 \put(0,4){\line(1,0){2}}
 \put(0,5){\line(1,0){2}}
  \put(0,0){\line(0,1){5}}
  \put(1,0){\line(0,1){3}}
  \put(2,3){\line(0,1){2}}
  \put(0.25,1.3){$s_0$}
  \put(0.25,2.35){$s_1$}
   \put(0.25,0.35){$s_2$}
 \put(0.25,3.35){$s_2$}
  \put(0.25,4.35){$s_0$}
    \put(1.25,4.35){$s_1$}
   \put(1.25,3.35){$s_0$}
  \end{picture}. 
  \end{center}
\end{ex}

\subsection{On the equivariant $K$-homology of $\Gr_{\Gad}$}
\label{SS:adjoint}
Let $G$ be a reductive group and $T$ the maximal torus. Let $(X,\Phi,X^\vee,\Phi^\vee)$ be the root data of $(G,T)$. 
Let $T_\ad$ be the maximal torus of $\Gad.$
We have $$
R(T_\ad)=\Z[Q]
$$
where 
$Q=\Z\Phi$ is the root lattice.
As the dual of the map $\phi: X^\vee\rightarrow P^\vee$ from \eqref{E:Xvee to Pvee},
there is a map $Q\rightarrow X$ of $\Z$-modules.  
So $R(T)=\Z[X]$ is an $R(T_\ad)$-algebra. 
Consider the diagram \eqref{E:fundamental as Dynkin auto}.
We regard $\ker(\overline{\phi})\subset X^\vee/Q^\vee$ as a subgroup of $\FG{G}$.
Note that there is a surjective map $G\rightarrow \Gad$ of algebraic groups.
\begin{prop}\label{P: Gr Gad} There is an injective homomorphism of $\Rep$-Hopf algebras
\begin{equation} \label{eq: K homology of adjoint group}
   K_*^T(\Gr_G)/J\rightarrow 
   R(T)\otimes_{R(T_\ad)} K_*^{T_\ad}(\Gr_{\Gad})
\end{equation}
where $J$ is the ideal generated by $\O_\fundel^G - 1$ for $\fundel\in \ker(\overline{\phi})$.
It is surjective if $\phi$ is surjective.
\end{prop}
\begin{proof}
The map $G\to \Gad$ induces a surjective map $\Gr_G\to\Gr_{\Gad}$. 
It is $T$-equivariant where $T$ acts on $\Gad$ via $T\mapsto T_\ad$. 
We denote the map $X^\vee(T)\to X^\vee(T_\ad)$ by $\la\mapsto \overline{\la}$.
There is an induced map
\begin{align*}\label{E:K G to Gad}
K^T_*(\Gr_G) &\overset{\Psi}\to R(T) \otimes_{R(T_\ad)} K^{T_\ad}(\Gr_{\Gad}) \\
i_{t_\la}^* &\mapsto i_{t_{\overline{\la}}}^*
\end{align*}
in the notation of Theorem \eqref{T:G iso}.

We will show that $\ker(\Psi)\supset J$.
For $\nu\in X^\vee$, by the definition of $\phi$ in \eqref{E:Xvee to Pvee},
$\nu\in\ker(\phi)$ if and only if $\nu$ pairs trivially with $Q$:
$\pair{\nu}{\alpha_i}=0$ for all $i\in I$. Call such $\nu$ \emph{central}.
Thus $\ker(\overline{\phi})$ consists of elements $\nu+Q^\vee$ where $\nu\in X^\vee$ is central. Let $\fundel\in \FG{G}$ correspond to $\nu+Q^\vee$ for $\nu\in X^\vee$ central, under the isomorphism \eqref{E:fund Xvee/Qvee}. We must show $\O_\fundel^G-1\in \ker(\Psi)$.
Considering the diagram \eqref{E:fundamental as Dynkin auto}, we see that
$\fundel$ induces the trivial automorphism of $\Diag$.
Therefore $\cl(\fundel)=1$ and $\fundel$ is a translation element, say, $\fundel=t_\gamma$ for $\gamma\in X^\vee$. $\gamma$ is central, for otherwise $t_\gamma$ has a descent, which is not possible since $\fundel=t_\gamma$ has length zero. 

By \eqref{E:fundamental basis} 
we have $\ell_\fundel = t_\gamma$. In the language of \eqref{E:G iso} we have
$i_{t_\gamma}^*=\O_{\fundel}^G$ and
$\Psi(\O_\fundel^G)=\O_{0}=1$. We deduce that $\O_\fundel^G-1\in \ker(\Psi)$.

For $\fundel\in \ker(\overline{\phi})$ and $w\in \hat{W}_G^0$, we have 
\begin{align*}
  \ell_{\fundel w} = \ell_\fundel \fundel(\ell_w) = \ell_\fundel \ell_w,
\end{align*}
where the second equality holds since $\cl(\fundel)=1$. This translates to the identity
\begin{equation}\label{eq: O pi factor}
    \O_{\fundel w}^G=\O_{\fundel}^G\O^G_{w}.
\end{equation}
Consider the composition 
\begin{equation}\label{eq: composition}
{ K_*^T(\Gr_G)}/J
\rightarrow K_*^T(\Gr_G)
/\mathrm{Ker}(\Psi)\rightarrow
R(T)\otimes_{R(T_\ad)} K_*^{T_\ad}(\Gr_{\Gad}).
\end{equation}
Let $\{\sigma_i\mid 1\le i\le n\}\subset \pi_1(G)$ be a set of representatives for $\pi_1(G)/\mathrm{Ker}(\overline{\phi}).$
Then in view of \eqref{eq: O pi factor}, $\{\O_{\sigma_i w}^G\;(1\le i\le n, \;w\in \hat{W}_G^0)\}$ is an $\Rep$-basis of   
${ K_*^T(\Gr_G)}/J$.
Since the image of this set 
is linearly independent over $\Rep$, the map \eqref{eq: composition} is injective and we have $J=\mathrm{Ker}(\overline{\phi}).$
Thus we have obtained an injective homomorphism \ref{eq: K homology of adjoint group} of $\Rep$-algebras.
It is clear the map preserves the Hopf-algebra structure. 
\end{proof}

\subsection{Affine Dynkin automorphisms fixing $0\in\Iaf$}\label{ssec:Affine Dynkin auto}
Let $\sigma\in \Aut(\Diag)$ be such that $\sigma(0)=0$.
By \eqref{E:pi on simple reflections} $\sigma$ defines a group automorphism of $\Waf$. Since
$\sigma$ is an automorphism for the finite Dynkin diagram, it
induces an automorphism on $Q$ by permuting simple roots
and an automorphism of $P^\vee$ by permuting fundamental 
coweights: 
$\sigma(\alpha_i)=\alpha_{\sigma(i)}$ 
and $\sigma(\om_i^\vee)=\om_{\sigma(i)}^\vee$ for $i\in I$.
These define actions of $\sigma$ on $X$, $X^\vee$,
and $R(T)$. The action on $P^\vee$ implies that
conjugation by $\sigma$ acts on $\FG{G}$ by a group automorphism.
Combining all this we obtain an involutive automorphism of 
$\Kx$ and $\Lx$. This action stabilizes $\Lx$ and the $k_w$ and $\ell_w$ bases.

\begin{prop} \label{P:auto homology} 
Let $\sigma\in\Aut(\Diag)$ be such that $\sigma(0)=0$. 
Then $\sigma$ induces an automorphism of $K^T_*(\Gr_G)$ 
given by 
\begin{align}
\label{E:sigma k}
\sigma(k_w) &= k_{\sigma(w)}, \\
\sigma(l_w) &= l_{\sigma(w)},
\end{align}
for all $w\in \Wxgr$.
\end{prop}
\begin{proof} For all $a\in\Kx$ and $b\in \Lx$ we have
\begin{align*}\label{E:sigma star}
  \sigma(a*b) = \sigma(a) * \sigma(b).
\end{align*}
We may assume $a = f t_\beta u$ and $b= g t_\gamma$ where
$f,g\in R(T)$, $u\in\Wfin$, and $\beta,\gamma\in X^\vee$. We have
\begin{align*}
  \sigma(a)*\sigma(b) &= \sigma(ft_\beta u)* \sigma(g t_\gamma) \\
  &= \sigma(f) t_{\sigma(\beta)} \sigma(u) * \sigma(g) t_{\sigma(\gamma)} \\
  &= \sigma(f u(g)) t_{\sigma(\beta)} t_{\sigma(u(\gamma))} \\
  &= \sigma(f u(g) t_\beta t_{u(\gamma)} \\
  &= \sigma(f t_\beta u * g t_\gamma) \\
  &= \sigma(a*b).
\end{align*}
We have $\sigma(D_w)=D_{\sigma(w)}$ for all $w\in \Wx$. Therefore
\begin{align*}
  \ell_{\sigma(w)} &= D_{\sigma(w)} * 1 \\
  &= \sigma(D_w)*\sigma(1) \\
  &= \sigma(D_w*1) \\
  &= \sigma(\ell_w).
\end{align*}
\end{proof}

\begin{rem} We list the nontrivial examples of $\sigma$.
For $A_{n-1}^{(1)}$, $\sigma(i)=-i$ mod $n$.
For $D_n^{(1)}$, $\sigma$ fixes all $i\in\Iaf$ except
$\sigma(n-1)=n$ and $\sigma(n)=n-1$. For $D_4^{(1)}$, $\sigma$ can be
any permutation of the three non-zero cominuscule nodes.
For $E_6^{(1)}$, $\sigma$ swaps the two branches not containing node $0$.
\end{rem}

\section{$K$-theoretic double $k$-Schur functions}\label{S:K k Schur}
In this section we introduce 
$\ksvar{k}{w}$ (resp. $ \tksvar{k}{w}$), the double $K$-theoretic (resp. closed) $k$-Schur functions. In \S \ref{SS:Pieri}, we 
study the double $K$-theoretic $k$-Schur functions
associated to the \emph{special} Schubert varieties.
In \S \ref{sec: k-conj}, we prove
a $k$-conjugation formula which is due to Morse for 
non-equivariant case. 
In \S\ref{SS: factorization}, we give the $k$-rectangle factorization of $\tksvar{k}{w}$. In \S \ref{SS: SL and PGL}, we explain
how the
Schubert bases for $G=\SL_n(\C), \PGL_n(\C)$ are also identified 
with the double $K$-theoretic $k$-Schur functions.

The $K$-theoretic formal group law $u \oplus v = u + v - uv$ will prove to be useful notationally.
The group law negation operation is $v\mapsto \ominus{v} = -v/(1- v)$. We write $u \ominus v = u \oplus (\ominus{v})=(u-v)/(1-v)$.

\subsection{Definition of $K$-theoretic double $k$-Schur functions}\label{S: Def K-k-Schur}
Let ${\Lambda}=\bigoplus_{i=0}^\infty \Lambda_i$ be the graded ring of
symmetric functions in 
the variables $y=(y_1,y_2,\ldots)$
with coefficients in $\Z$ (\cite{Mac}).
 Let $\hat{\Lambda}$ be the 
completed ring of symmetric functions in $y$, which consists of formal sums $f_0+f_1+\cdots$ where $f_i\in \Lambda_i.$ For any commutative ring $A$,
let $\hL{A}:=A\otimes_\Z\hat{\Lambda}.$

For $1\le i\le n$, define
\begin{equation}
b_i=1-e^{-a_i} = c(-a_i)\in \RepGL.
\end{equation}
The following power series in $\hL{\Rep}$ plays a fundamental 
role in the construction of the $K$-theoretic double $k$-Schur functions:
\begin{align*}
 \Omega(b_i|y)&=\prod_{j=1}^{\infty}(1-b_i y_j)^{-1} \in \hL{\RepGL}.
\end{align*}

\begin{prop}\label{P: KGL acts on hat Lambda}
$\hL{\RepGL}$ is a $\KGL$-module, extending the
action of $\WGL$ such that
$\Wfin=S_n$ acts on the
coefficient ring $\RepGL$ by permuting $(b_1,\dotsc,b_n)$ 
and the generators $t_{\pm \eps_i}$ of $X^\vee\subset \Wx$ act by multiplication by 
$\Omega(b_i|y)^{\pm1}$ for $1\le i\le n$. In particular
\begin{align*}
  s_0 (f(y|b))&=\frac{\Omega(b_1|y)}{\Omega(b_n|y)} f(y|s_\theta(b)),
  \\
  \sh( f(y|b))&=\Omega(b_1|y)f(y|\rot(b)),\label{eq:sh action}
\end{align*}
for $f(y|b)\in \hL{\RepGL}$. Here $\rot=u_1=s_1s_2\dotsm s_{n-1}$, so that
$\rot(b_i)=b_{i+1}$ for $1\le i\le n-1$ with the convention $b_{n+1}=b_1$.
\end{prop}
\begin{proof} 
For $u\in\Wfin$ and $1\le i\le n$ we have
$u t_{\pm \eps_i} u^{-1} = t_{\pm \eps_{u(i)}}$.
Write $f(y|b)\in \hLa_\Rep$ for a symmetric series in $y$ with coefficients depending on 
$b=(b_1,\dotsc,b_n)$. We have
\begin{align*}
  u t_{\pm \eps_i} u^{-1} \cdot f(y|b) &= u\cdot \Omega(b_i|y)^{\pm1} f(y|u^{-1}(b)) \\
  &= \Omega(b_{u(i)}|y)^{\pm1} f(y|b) \\
  &= t_{\pm \eps_{u(i)}}\cdot f(y|b).
\end{align*}
This well-defines an action of $\KGL^\reg$ on $\hLa_{\RepGL^\reg}$.
To show this restricts to an action of $\KGL$ on $\hLa_{\RepGL}$, 
it must be shown that $T_i(\hL{\RepGL})\subset \hL{\RepGL}$, 
This holds for $i\in I$ since the action of $T_i$ is only on the coefficient ring $\RepGL$.
For $i=0$ we have
\[
c(\al_0)\Omega(b_n|y)T_0(f)=
\Omega(b_1|y)s_\theta(f)-\Omega(b_n|y)f=
(s_\theta-1)\left(
\Omega(b_n|y)f
\right).
\]
Because 
$(s_\theta-1)
(\Omega(b_n|y)f)$ is divisible by $b_1-b_n=c(\al_0)(1-b_n),$
and $1-b_n$ and $\Omega(b_n|y)$ are invertible, we have $T_0(f)\in \hL{\RepGL}.$
\end{proof}

$\hL{\Rep}$ has the structure of an $\Rep$-Hopf algebra coming from 
symmetric functions: 
the coproduct $\Delta$,
the antipode $S$, and 
counit $\epsilon$ are given by 
$$
\Delta(h_r)=\sum_{i=0}^{r}h_i\otimes h_{r-i},\quad
S(h_r) = -e_r, \quad\text{and} \quad\epsilon(h_r)=0\quad\text{for $r\in\Z_{>0}$,}
$$
where $h_r$ (resp. $e_r$) denotes the $r$-th complete (resp. elementary) symmetric function.
It is straightforward to verify
\begin{equation}
\label{E: Hopf structure on hat Lambda}
\Delta(\Omega(b_i|y))=\Omega(b_i|y)
\otimes\Omega(b_i|y),\quad
S(\Omega(b_i|y))=
\Omega(b_i|y)^{-1},\quad
\epsilon(\Omega(b_i|y))=1,
\end{equation}
for $1\le i\le n.$

\begin{defn}
Let $\hLa_{(n)} = \KGL\cdot 1\subset \hL{\Rep}$. For $w\in\WGLgr$, the
\emph{$K$-theoretic double $k$-Schur function} $g^{(k)}_w(y|b)$ and
\emph{closed $K$-theoretic double $k$-Schur function} $\tilde{g}^{(k)}_w(y|b)$
are the elements of $\hLa_{(n)}$ defined by
\begin{align}
g^{(k)}_w(y|b) &= T_w(1), \\
\tilde{g}^{(k)}_w(y|b) &= D_w(1).
\end{align}
\end{defn}

Using \eqref{E:T basis} and the fact that $T_u(1)=0$ for $u\in \Wfin$, we deduce that
$\hLa_{(n)}$ is spanned over $\Rep$ by $\{g^{(k)}_w(y|b)\mid w\in \WGLgr\}$, which connects with the definition of $\hLa_{(n)}$ given in the introduction. 
We will show that this spanning set is a basis.

\begin{lem}\label{L:Jac}
The elements $\Omega(b_1|y),\ldots,\Omega(b_n|y)$ are algebraically 
independent over $\mathrm{Frac}(\Rep)$.
\end{lem}
\begin{proof}
It suffices to show the algebraic independence for $\Omega(b_i|y)$ specialized to the variables $y_1,\ldots,y_n.$
It is straightforward to 
show that $$
\det(\partial \Omega(b_i|y)/\partial y_j)_{1\le i,j\le n}
=\prod_{i=1}^nb_i\Omega(b_i|y)\times
\prod_{i<j}(b_i-b_j)
(y_i-y_j).$$
Hence the proposition holds. 
\end{proof}

\begin{proof}[Proof of Theorem \ref{thm:main}]
By Theorem \ref{T:G iso} we may replace $K^T_*(\Gr_{\GL_n})$ by $\LGL$
and $\cI_w$ and $\O_w$ by $k_w $ and $\ell_w$ respectively.
By \eqref{E:L-reg} there is a $\Rep^\reg$-algebra homomorphism
\begin{align}
\LGL^\reg =\Rep^\reg[X^\vee]&\overset{\psi^\reg}\longrightarrow
\hat{\Lambda}_{\Rep^\reg},
\label{E:psi-reg} \\
\label{E:psi}
t_{\pm \eps_i}&\mapsto\Omega(b_i|y)^{\pm 1}\quad \text{for $1\le i\le n$}.
\end{align}
It is injective by Lemma  \ref{L:Jac}.
It is straightforward to 
verify that $\psi^\reg$ is $\WGL$-equivariant and therefore $\KGL$-equivariant.
Let $\psi:\LGL \to \hLa_{(n)}$ be the restriction of $\psi^\reg$. We have
\begin{align}\label{E:psi O}
\psi(k_w)=\psi(T_w*1)=T_w(1)=\tksvar{k}{w}.
\end{align}
Similarly $\psi(\ell_w)=\tksvar{k}{w}$.
We see that the image of $\psi$ is $\hLa_{(n)}$. 

It is clear that $\psi$ preserves the Hopf structure maps given by \eqref{E: Hopf structure on hat Lambda}.
\end{proof}

\begin{rem} \label{R:Omegas are GL classes}
For $1\le i\le n$, $\Omega(b_i|y)^{\pm1}\in \hLa_{(n)}$
by \eqref{E:psi} since $t_{\pm\epsilon_i}\in \LGL$.
\end{rem}

\begin{rem} \label{R:Wfin action}
Under the isomorphism \eqref{E:G iso},
the action of $\KGL$ on $K^T_*(\Gr_\GL)$ 
corresponds to the $*$-action of $\KGL$ on $\LGL$.
Restricted to $\Wfin\subset \KGL$, 
the corresponding $\Wfin$ action on $\hLGL$ permutes the parameters $b_i$.
\end{rem}

\begin{prop}\label{P: double K-k for sh}
For $i\in \Z$, we have
\[g_{\sh^i}^{(k)}(y|b)=
\tilde{g}_{\sh^i}^{(k)}(y|b)
=\begin{cases}
    \Omega(b_{1}|y)
    \Omega(b_2|b)
\cdots\Omega(b_{i}|y) & \text{if $i>0$}\\
 1 & \text{if $i=0$}\\
\Omega(b_{i+1}|y)^{-1}\cdots\Omega(b_{-1}|y)^{-1}\Omega(b_{0}|y)^{-1} &\text{if $i<0$}
\end{cases},
\]
where we put $b_{i+nl}=b_i$ for $l\in \Z.$
For $i\in \Z$ and $w\in \WSLgr$,
we have
\begin{align}
\ksvar{k}{\sh^i\, w}\label{eq:sh acts on K-k-Schur}
&=\ksvar{k}{\sh^i}\, \ks{k}{w}(y|\rot^{i}(b)),\\
\tilde{g}_{\sh^i\, w}^{(k)}(y|b)
&=\ksvar{k}{\sh^i}\,\tks{k}{w}(y|\rot^{i}(b)).
\end{align}
\end{prop}

\begin{rem}\label{R: O sh n}
Since $\rot^n=1$, we have from Proposition \ref{P: double K-k for sh}
$$
\tksvar{k}{\sh^n w}=\tksvar{k}{\sh^n}\tksvar{k}{w}
=\Omega(b_1|y)\cdots \Omega(b_n|y)\tksvar{k}{w}.
$$
This identity corresponds to $\O_{\sh^n w}=
\O_{\sh^n}\O_{w}$ for $w\in \WGLgr.$
\end{rem}

\begin{cor}
For $i\in \{0,1,\ldots,n-1\}$, $w\in \WGLgr$, and $\fundel\in\FG{\GL_n}$
\begin{align*}
T_i \,\ksvar{k}{w}&=
\begin{cases}
\ksvar{k}{s_iw}& (s_iw>w \;\text{and $s_iw\in \WGLgr$})\\
0 & (s_iw>w \;\text{and $s_iw\notin \WGLgr$})\\
-\ksvar{k}{w}& (s_iw<w)
\end{cases},
\\
D_i \,\tksvar{k}{w}&=
\begin{cases}
\tksvar{k}{s_iw} & (s_iw>w \;\text{and $s_iw\in \WGLgr$} )\\
\tksvar{k}{w}& (\text{otherwise})
\end{cases},
\\
\fundel \, \ksvar{k}{w} &= \ksvar{k}{\fundel w}, \\
\fundel \, \tksvar{k}{w} &= \tksvar{k}{\fundel w}.
\end{align*}
\end{cor}

\begin{proof}
This follows from Proposition \ref{P:T on k} and Corollary \ref{C:TD on kl}.
\end{proof}

\subsection{Special Schubert classes}\label{SS:Pieri}
For $1\le i\le k$ define the affine Grassmann elements $\rho_i$ and $\rho'_i$ by
\begin{align}
\label{E:single row}
\rho_i&=s_{i-1}s_{i-2}\dotsm s_1s_0, \\
\label{E:single column}
\rho_i'&=s_{n-i+1} \dotsm s_{n-1}s_0.
\end{align}
The Schubert varieties associated to 
$\rho_i$ and $\rho_i'$ are called \emph{special}
Schubert varieties and the corresponding
classes are called the Pieri classes.
We study the double $K$-$k$-Schur functions 
associated to the Pieri classes.

\subsubsection{Plethystic notation}
In this section some plethystic notation is used for symmetric series or polynomials. For sets of variables $A$ and $B$ let $h_k[A-B]$ be defined by the generating series
\begin{align*}
  \sum_{r\ge0} h_r[A-B] t^r = \dfrac{\prod_{b\in B} (1-tb)}{\prod_{a\in A} (1-ta)}
\end{align*}
Since the ring $\La$ of symmetric functions is a polynomial ring with generators $h_k$ for $k\ge1$, there is a ring homomorphism denoted by $f\mapsto f[A-B]$ for $f\in \La$ and defined by
$h_k\mapsto h_r[A-B]$ for $r\ge1$, from $\La$ to polynomials/series that are symmetric in both $A$ and $B$. We have
\begin{align*}
s_\la[A-B] &= \sum_{\mu\subset\la} s_\mu[A] s_{\lambda/\mu}[-B] \\
&= \sum_{\mu\subset\la} (-1)^{|\la|-|\mu|} s_\mu[A] s_{\lambda'/\mu'}[B]
\end{align*}
where $s_{\la/\mu}$ is the skew Schur function and $\la'$ is the conjugate or transpose partition of $\la$.

Let $[i,j]$ denote the interval of integers $\{i,i+1,\dotsc,j\}$.
For an interval of integers $I$ and $p\in I$ let $I\setminus p$ be the
interval $I$ with the element $p$ removed. Finally, for a set of integers $I$
let $b_I$ be the set of variables $b_i$ for $i\in I$.

\subsubsection{Schur expansion of the special double $K$-$k$-Schur functions}

\begin{prop}\label{prop:single-row} For $1\le i\le k$,
\begin{align}
\ksvar{k}{\rho_i}&=
e^{-(a_1+\cdots+a_{i-1}
+a_n)}\sum_{p\ge1,q\ge0}(-1)^q b_n^q 
 h_{p-i}(b_1,\dots,b_i) s_{(p,1^q)}(y),\label{eq:singlerow}\\ 
\label{eq:singlecolumn}
\ksvar{k}{\rho'_i}&=
e^{-a_{n-i+1}}\sum_{p\ge1,q\ge0}(-1)^{q+i-1} 
b_1^{p-1} \times \\ \notag
&\qquad\left[\sum_{r=0}^{i-1}(-1)^{i-1-r}\binom{i-1}{r} h_{q-r}(b_{n-i+1},\dots,b_n)\right]
s_{(p,1^q)}(y).
\end{align}
and in particular we have
\begin{align}
\label{E:single row non equiv}
\ks{k}{\rho_i}(y|0) &=h_i(y), \\
\label{E:single column non equiv}
\ks{k}{\rho'_i}(y|0) &= \sum_{r=0}^{i-1}  \binom{i-1}{r} e_{r+1}(y).
\end{align}
\end{prop}

\begin{proof}
The proof proceeds by induction on $i$.
In the case of $i=1$, we know that \ref{eq:singlerow} is true by the following.
\begin{align*}
  c(\alpha_0) T_0(1) &= -1 + \Omega[(b_1-b_n)y] \\
  &= \sum_{\lambda\ne\emptyset} s_\lambda[b_1-b_n] s_\lambda[y] \\
  &= \sum_{p\ge1, q\ge0} s_{(p,1^q)}[b_1-b_n] s_{(p,1^q)}[y] \\
  &= \sum_{p\ge1, q\ge0} (b_1^p(-b_n)^q +b_1^{p-1}(-b_n)^{q+1}) s_{(p,1^q)}[y] \\
  &= (b_1-b_n) \sum_{p\ge1,q\ge0} b_1^{p-1} (-b_n)^q s_{(p,1^q)}[y].
\end{align*}
This implies \eqref{eq:singlerow} for $i=1$.

Next, assume that \eqref{eq:singlerow} is true when $i-1>0$. In this case, the following
calculation shows that \eqref{eq:singlerow} is also true for $i$, using the identity
\begin{align*}
(b-a) h_r[A] = h_{r+1}[A-a]-h_{r+1}[A-b]
\end{align*}
for any variable set $A$ and variables $a$ and $b$, with $A=[i]$, $a=b_{i-1}$ and $b=b_i$.
We have
    \begin{align*}
        &c(\al_{i-1})
        \ksvar{k}{\rho_{i}} \\
        &=s_{i-1}(\ksvar{k}{\rho_{i-1}})-\ksvar{k}{\rho_{i-1}}\\
        &=e^{-(a_1+\cdots+a_{i-2}+a_n)}\sum_{p,q\geq0}(-1)^q b_n^q \left( h_{p-i+1}(b_{[i]\setminus \{i-1\}}) -h_{p-i+1}(b_{[i]\setminus\{ i\}}) \right)s_{(p,1^q)}(y)\\
&=
e^{-(a_1+\cdots+a_{i-2}+a_n)}(b_i-b_{i-1})\sum_{p,q\geq0}(-1)^q b_n^q h_{p-i}(b_1,\dots,b_{i}) s_{(p,1^q)}(y).
\end{align*}
Thus we have
\begin{align*}
\ksvar{k}{\rho_i}&=
e^{-(a_1+\cdots+a_{i-1}+a_n)}\sum_{p,q\geq0}(-1)^q b_n^q h_{p-i}(b_1,\dots,b_{i}) s_{(p,1^q)}(y)
    \end{align*}
proving \eqref{eq:singlerow}.

Next, \eqref{eq:singlecolumn} is shown by induction on $i$ as in the proof of \ref{eq:singlerow}.
For $i=1$, $\rho'_1=\rho_1$ and the formulas \eqref{eq:singlerow} and \eqref{eq:singlecolumn} agree for $i=1$. It follows that \eqref{eq:singlecolumn} holds for $i=1$.

We assume that \eqref{eq:singlecolumn} is true when $i-1>0$. In this case, the following
calculation shows that \eqref{eq:singlecolumn} is also true for $i$.
We use the identity
\begin{align*}
    (1-b) h_p[A-a] - (1-a) h_p[A-b] = (b-a) (h_{p-1}[A]-h_p[A])
\end{align*}
with $b=b_{n-i+1}$, $a=b_{n-i+2}$, and $A=b_{[n-i+1,n]}$.
We have
\begin{align*}
&c(\al_{n-i+1})\ksvar{k}{\rho'_{i}} \\
&=s_{n-i+1}\left( \ks{k}{\rho'_{i-1}}\right)-
\ks{k}{\rho'_{i-1}}\\
&=\sum_{p\ge1,q\ge0}(-1)^{q+i-2} 
b_1^{p-1} s_{(p,1^q)}(y)\left[\sum_{r=0}^{i-2}(-1)^{i-2-r}\binom{i-2}{r} \times \right. \\
&\qquad\Bigl((1-b_{n-i+1}) h_{q-r}(b_{[n-i+1,n]\setminus \{n-i+2\}})-(1-b_{n-i+2})h_{q-r}(b_{[n-i+1,n]\setminus\{ n-i+1\}}
\Bigr) \Bigl.\Bigr]
\\
&=(b_{n-i+1}-b_{n-i+2})\sum_{p\ge1,q\ge0}(-1)^{q+i-2} 
b_1^{p-1}s_{(p,1^q)}(y) \times \\
&\qquad\sum_{r=0}^{i-2}(-1)^{i-2-r}\binom{i-2}{r}
\Bigl[(h_{q-r-1}(b_{[n-i+1,n]}) -h_{q-r}(b_{[n-i+1,n]})) \Bigr]. \\
\end{align*}
Thus we have
\begin{align*}
\ksvar{k}{\rho'_{i}} &= e^{-a_{n-i+1}} \sum_{p\ge1,q\ge0}(-1)^{q+i-1} 
b_1^{p-1}s_{(p,1^q)}(y) \times \\
&\qquad\sum_{r=0}^{i-1}(-1)^{i-1-r}\left( \binom{i-2}{r} + \binom{i-2}{r-1}\right)
h _{q-r}(b_{[n-i+1,n]}))  \\
&= e^{-a_{n-i+1}} \sum_{p\ge1,q\ge0}(-1)^{q+i-1} 
b_1^{p-1}s_{(p,1^q)}(y) \times \\
&\qquad\sum_{r=0}^{i-1}(-1)^{i-1-r}\binom{i-1}{r} h _{q-r}(b_{[n-i+1,n]})) 
\end{align*}
proving \eqref{eq:singlecolumn}.
\end{proof}

\subsubsection{Pieri classes and the translation elements}
\begin{prop}\label{prop:Om-Pieri}
For $1\le i\le k$ we have
\begin{align}
\frac{\Omega(b_i|y)}{\Omega(b_n|y)}
&=1+c(a_n-a_i)\sum_{j=1}^{i}
\prod_{l=1}^{j-1}c(a_l-a_i)
g_{\rho_j}^{(k)}(y|b),\label{eq:Om-Pieri}\\
\frac{\Omega(b_1|y)}{\Omega(b_i|y)}
&=1+c(a_i-a_1)\sum_{j=1}^{n-i+1}
\prod_{l=i+1}^{n}c(a_i-a_l)
g_{\rho_j'}^{(k)}(y|b).\label{eq:Om-Pieri-col}
\end{align}
\end{prop}
\begin{proof} The proof is by induction on $i$.
For $i=1$ we have
    \begin{align*}
    \ksvar{k}{\rho_1} &= T_0(1) = c(\al_0)^{-1} \left(\frac{\Omega(b_1|y)}{\Omega(b_n|y)}-1\right) 
     \end{align*}
     which implies \eqref{eq:Om-Pieri} for $i=1$ since $\alpha_0=a_n-a_1$. 
     
    Next, assume that the claim is true for $i-1>0$. In this case, the following calculation shows that the claim is also true for $i$.
    \begin{align*}
        \frac{\Omega(b_i|y)}{\Omega(b_n|y)}
        &=s_{i-1}\left( \frac{\Omega(b_{i-1}|y)}{\Omega(b_n|y)}\right)
          =s_{i-1}
          \left( 1+c(a_n-a_{i-1})\sum_{j=1}^{i-1}\prod_{l=1}^{j-1}c(a_l-a_{i-1})g_{\rho_j}^{(k)}(y|b)
          \right)\\
          &=1+c(s_{i-1}(a_n-a_{i-1}))\sum_{j=1}^{i-1}\prod_{l=1}^{j-1}c(s_{i-1}(a_l-a_{i-1}))s_{i-1}(g_{\rho_j}^{(k)}(y|b))\\
          &=1+c(a_n-a_{i})\sum_{j=1}^{i-1}\prod_{l=1}^{j-1}c(a_l-a_{i})s_{i-1}(g_{\rho_j}^{(k)}(y|b))\\
          &=1+c(a_n-a_i)
          \left(
          \sum_{j=1}^{i-2}\prod_{l=1}^{j-1}c(a_l-a_{i})g_{\rho_j}^{(k)}(y|b)+
          \prod_{l=1}^{i-2}c(a_l-a_{i})s_{i-1}(g_{\rho_{i-1}}^{(k)}(y|b))
          \right)\\
          &=1+c(a_n-a_{i})\\
          &\left(
          \sum_{j=1}^{i-2}\prod_{l=1}^{j-1}c(a_l-a_i)g_{\rho_j}^{(k)}(y|b)+
          \prod_{l=1}^{i-2}c(a_l-a_i)(g_{\rho_{i-1}}^{(k)}(y|b)+c(a_{i-1}-a_i)g_{\rho_{i}}^{(k)}(y|b))
          \right)\\
          &=1+c(a_n-a_i)\sum_{j=1}^{i}\prod_{l=1}^{j-1}c(a_l-a_i)g_{\rho_j}^{(k)}(y|b).
    \end{align*}
    \eqref{eq:Om-Pieri-col} is proved similarly.
\end{proof}

\subsubsection{Forgetting equivariance}
The goal of this subsection is Corollary \ref{C:forget to K k schur}.

Let $\phi_0:\Rep\to\Z$ be the ring homomorphism $\phi_0(e^\la)=1$ for all $\la\in X$.
For any commutative $\Rep$-algebra $A$, there is a ring homomorphism $A\to \Z \otimes_{\Rep} A$ given by $a\mapsto 1\otimes a$ using $\phi_0$ for the action of $\Rep$ on $\Z$.
Define $\phi:K^T_*(\Gr_G) \to \Z \otimes_{\Rep} K^T_*(\Gr_G) \cong K_*(\Gr_G)$.

For $u,v,w\in \Wxgr$ define the structure constants $\cc_{uv}^w\in \Rep$ by
\begin{align}\label{E:T constants}
  \cI_u \cI_v &= \sum_{w\in \Wxgr} \cc_{uv}^w \cI_w.
\end{align}
Applying the surjective ring homomorphism $\phi$ we have 
\begin{align}\label{E:Z constants}
  \phi(\cI_u)\phi(\cI_v) &= \sum_{w\in \Wxgr} \phi_0(\cc_{uv}^w) \phi(\cI_w).
\end{align}

We now consider $\Gr_{\SL_n}$ but use scalars in the representation ring of $T\subset \GL_n(\C)$.
In this case we rewrite $\cc_{uv}^w$ as $\cc_{uv}^w(b)$ and
$\phi_0(\cc_{uv}^w)$ as $\cc_{uv}^w(0)$ since $\phi_0(b_i)=0$ for all $1\le i\le n$.

Let $\psi:K^T_*(\Gr_{\GL_n})\to \hLa_{(n)}$ be the isomorphism in Theorem \ref{thm:main}
and let $\psi_0:K_*(\Gr_{\SL_n})\cong \La_{(n)}$ denote the isomorphism in \cite[Thm. 7.17]{LSS:K}.
Let $\phi_{\hL{}}$ be the map
$\hL{\Rep} \to \Z \otimes_{\Rep} \hL{\Rep} \cong \hL{}$.
We have a diagram of ring homomorphisms
\[
\begin{tikzcd}
K^T_*(\Gr_{\GL_n}) \arrow[d,swap,"\psi"] &K^T_*(\Gr_{\SL_n})  \arrow[r,"\phi"] \arrow[l,swap,"\supseteq"] & K_*(\Gr_{\SL_n}) \arrow[d,"\psi_0"] \\
\hLa_{(n)} \arrow[r,"\phi_{\hL{}}"]& \hL{} & \La_{(n)}\arrow[l,swap,"\supseteq"]
\end{tikzcd}
\]

\begin{cor}\label{C:forget to K k schur}
For $w\in \hat{W}_{\SL_n}^0$, $\ks{k}{w}(y|0)$ is the $K$-theoretic $k$-Schur function $g_w^{(k)}(y)$ of \cite{LSS:K}.
\end{cor}
\begin{proof} By \cite[Thm. 7.17]{LSS:K} and Remark \ref{R:KGr},
$\psi_0(\phi(\cI_w))=g^{(k)}_w(y)$.
We have $\phi_{\hLa{}}(\psi(\cI_w))=g^{(k)}_w(y|0)$ by Theorem \ref{thm:main}
Applying the ring homomorphism $\psi_0\circ\phi$ to \eqref{E:T constants},
we obtain
\begin{align*}
    g_u^{(k)}(y) g_v^{(k)}(y) = \sum_{w\in \hat{W}_{\SL_n}^0} \cc_{uv}^w(0) g_w^{(k)}(y).
\end{align*}
Applying $\phi_{\hLa}\circ \psi$ to \eqref{E:T constants}, 
by Theorem \ref{thm:main} $\cI_w\mapsto g_w^{(k)}(y|0)$ and we have
\begin{align}
\label{E:g0 structure constants} g^{(k)}_u(y|0) g^{(k)}_v(y|0) &= \sum_{w\in \hat{W}_{\SL_n}^0} \cc^w_{uv}(0) g^{(k)}_w(y|0) .
\end{align}
By \eqref{E:single row non equiv} and \cite[Cor. 7.18]{LSS:K} we have
\begin{align}
g^{(k)}_{\rho_i}(y|0) = h_i(y) = g^{(k)}_{\rho_i}(y)\qquad\text{for $1\le i\le n-1$.}
\end{align}
Since the $h_i(y)$ for $1\le i\le n-1$ generate $\La_{(n)}$ and multiply against 
the $g^{(k)}_w(y|0)$ and $g^{(k)}_w(y)$ the same way, we deduce that
$g^{(k)}_w(y|0)=g^{(k)}_w(y)$ as required.
\end{proof}

\begin{rem} It is not directly obvious that the specializations $g^{(k)}_w(y|0)$ have bounded degree, but this follows from the above proof. 
\end{rem}

\subsection{$K$-theoretic $k$-conjugation on $\hL{\Rep}$}
\label{sec: k-conj}

A symmetry property of $k$-Schur functions called $k$-\emph{conjugation} was discovered by Lapointe and Morse in \cite{LM3}. Morse \cite[Theorem 61]{Mor} 
generalized the property for the $K$-theoretic $k$-Schur functions.
We apply the results in \S \ref{ssec:Affine Dynkin auto} to show the equivariant version of this property. 

\subsubsection{Involution $\hat{\iota}$}
Let $\iota$ be the automorphism of type $A_{n-1}^{(1)}$ given by $\iota(i)=-i$ mod $n\Z$ for $i\in \Iaf$. Denote the induced group automorphism of $\WSL$ by $w\mapsto w^*$. Then $(s_i)^*=s_{i^*}=s_{-i}$ with indices taken mod $n\Z$. This automorphism extends to $\WGL$ via $\sh^*=\sh^{-1}$. 
The induced action on the root lattice $Q$ is given by $\iota(\alpha_i)=\alpha_{n-i}$ for $i\in I$. 
The action extends to $X(\TGL)$ so that $\iota(a_i)=-a_{n+1-i}$ for $1\le i\le n$, which is equivalent to
$\iota(b_i) = \ominus b_{n+1-i}$. 
There is an induced involution on $K_*^T(\Gr_{\GL_n})$ denoted by the same symbol $\iota$ (see \S \ref{ssec:Affine Dynkin auto} ).

For $w\in \WGLgr$, we have from Proposition \ref{P:auto homology}
    \begin{equation}\label{eq: sigma on O}
    \iota(k_w)=k_{w^*}, \quad 
    \iota(\ell_w)=\ell_{w^*}.
    \end{equation}
Let us translate the involution $\iota$ on $\hL{\Rep}.$
Let $\kconj$ be a ring automorphism of $\Lambda$ defined in \cite[\S 8]{Mor} by 
$$\kconj(h_l(y))=\sum_{r=0}^{i-1}  \binom{i-1}{r} e_{r+1}(y)$$ for $l\ge 1.$
Note that $\kconj(h_l(y))=g_{\rho_i'}^{(k)}(y|0)$ by \eqref{E:single column non equiv}.
Equivalently, we have 
\begin{equation}\label{eq: K-k-conj}
\left(\sum_{l=0}^\infty u^l h_l(y)\right)\left(\sum_{l=0}^\infty \left(\ominus u\right)^l\kconj(h_l(y))\right)=1,
\end{equation}
where $u$ is an indeterminate.
It is known that $\kconj^2=\mathrm{id}$ (see \cite[Remark 59]{Mor} for a proof by A. Garsia). 

Let $\hat{\iota}$ be the ring automorphism of 
$\hL{\Rep}$ such that 
\begin{align*}
\hat{\iota}(h_l(y))&=\kconj h_l(y)\quad (l\ge 1),\\
    \hat{\iota}(b_i)&=\ominus b_{n+1-i}\quad (1\le i\le n).
\end{align*}
We have
$\hat{\iota}^2=\mathrm{id}$.
\begin{lem}\label{L: sigma on Omega}
We have
  \begin{equation} \label{E: sigma on Omega}
      \hat{\iota}(
\Omega(b_i|y))= \Omega(b_{n+1-i}|y)^{-1}\quad (1\le i\le n).
  \end{equation}
\end{lem}
\begin{proof}
We  have
 \begin{align*}
\hat{\iota}(\Omega(b_i|y))&=\hat{\iota}\left(\sum_{l=0}^\infty b_i^l h_l(y)\right)
= \sum_{l=0}^\infty (\ominus b_{n+1-i})^l\kconj (h_l(y))\\
&=\left(
\sum_{l=0}^\infty b_{n+1-i}^l h_l(y)\right)^{-1}
=\Omega(b_{n+i-i}|y)^{-1},
\end{align*}   
where we used \eqref{eq: K-k-conj} in the third equality.  
\end{proof}

\begin{prop} By the isomorphism $K_*^T(\Gr_{\GL_n})\cong\hLGL $,
the automorphism $\iota$
on $K_*^T(\Gr_{\GL_n})$ is transported 
to $\hat{\iota}$ on $\hLGL.$ Therefore we have
\begin{equation}\label{E: sigma on K-k-Schur}
    \hat{\iota}(\tksvar{k}{w})=\tks{k}{w^*}(y|b),\quad
     \hat{\iota}(\ksvar{k}{w})=\ks{k}{w^*}(y|b).
    \end{equation}
\end{prop}
\begin{proof} The automorphism $\iota$ on $K_*^T(\Gr_{\GL_n})\cong \LGL$ is determined by 
$\iota(t_{\eps_i})=t_{-\eps_{n+1-i}}$ 
or equivalently by $\iota(b_i)=\ominus b_{n+1-i}$
for $1\le i\le n.$ Hence by Lemma \ref{L: sigma on Omega}, the first statement holds. 
The second statement follows from \eqref{eq: sigma on O}.
\end{proof}

\begin{rem}
    It is not difficult to 
    show $\hat{\iota} \circ T_i
=T_{n-i} \circ \hat{\iota},\;\hat{\iota} \circ D_i
=D_{n-i} \circ \hat{\iota}$ directly. Then 
\eqref{E: sigma on K-k-Schur} follows immediately.
\end{rem}
\begin{rem} For $i\in \Z$, we have
  \begin{equation}
    \hat{\iota}(\tks{k}{w}(y|\rot^i b))=\tks{k}{w^*}(y|\rot^{-i}b),\quad
     \hat{\iota}(\ks{k}{w}(y|\rot^i b))=\ks{k}{w^*}(y|\rot^{-i}b).
    \end{equation}  
\end{rem}

\subsubsection{$\kconj$ on power sums}
In \cite{LLS:backstableK}, $\kconj$ was given in terms of power sum symmetric functions $p_j(y)$
(see \cite[(4.2)]{LLS:backstableK}). 
For the reader's convenience, we record a proof of the coincidence of two descriptions.  
\begin{prop}
\begin{align}\label{E:K transpose}
  \kconj(p_j(y)) = (-1)^{j+1} \sum_{r\ge0} \binom{r+j-1}{j-1} p_{j+r}(y)\qquad\text{for $j\ge1$.}
\end{align}
\end{prop}
\begin{proof}Newton's formula reads
\begin{align*}
\sum_{l=0}^\infty u^l h_l
=\exp\left(
\sum_{l=1}^\infty \frac{p_l}{l}u^l
\right),\quad
\sum_{l=0}^\infty (-u)^l e_l
=\exp\left(-
\sum_{l=1}^\infty \frac{p_l}{l}u^l
\right).
\end{align*}
From \eqref{eq: K-k-conj}, we have
\begin{align*}
\sum_{l=0}^\infty (-u)^l e_l
=\sum_{l=0}^\infty \left(\ominus u\right)^l\kconj(h_l)
=\exp\left(
\sum_{l=1}^\infty \frac{\kconj(p_l)}{l}(\ominus u)^l
\right)
\end{align*}
Therefore
\begin{align*}\exp\left(-
\sum_{l=1}^\infty \frac{p_l}{l}u^l
\right)
=\exp\left(
\sum_{l=1}^\infty \frac{\kconj (p_l)}{l}(\ominus u)^l
\right),
\end{align*}
or equivalently, by substitution $u\mapsto \ominus u$, 
\begin{align*}\exp\left(-
\sum_{l=1}^\infty \frac{p_l}{l}(\ominus u)^l
\right)
=\exp\left(
\sum_{l=1}^\infty \frac{\kconj (p_l)}{l}u^l
\right).
\end{align*}
It follows that 
\begin{align*}
-
\sum_{l=1}^\infty \frac{p_l}{l}(\ominus u)^l
=
\sum_{l=1}^\infty \frac{\kconj(p_l)}{l}u^l,
\end{align*}
which is equivalent to
\eqref{E:K transpose}.
\end{proof}

\subsection{$k$-rectangle factorization for $\tksvar{k}{w}$}
\label{SS: factorization}
Recall that we set $k=n-1$.
Let $\Par^k$ be the set of partitions $\lambda$ such that $\la_1\le k.$
There is a bijection (\cite{LLMSSZ})
\begin{align*}
\Par^k \to \WSLgr, \quad
\la \mapsto x_\la = {\pi_1^{-l+1}}(\rho_{\la_l}) \cdots {\pi_1^{-1}}(\rho_{\la_2}) \rho_{\la_1}
\end{align*}
for $\la=(\la_1,\ldots,\la_l)$.

Recall that in 
$\tilde{W}_{\PGL_n}^0$, we have
\begin{equation*}
t_{-\om_i^\vee}=\pi_i^{-1}\kappa_i,\quad
\kappa_i=x_{R_i}.
\end{equation*}

For any partitions $\la,\mu$, let
$\la\cup\mu$ denote
the partition obtained by concatenating $\la,\mu$ as sequences of non-negative integers and reordering in weakly decreasing order.
We obtain an equivariant analogue of the nonequivariant $K$-homology $k$-rectangle factorization due to Takigiku \cite{Tak}.
\begin{thm}[\cite{IINY}] \label{C:factor kappa in function}
For any $\la\in\Par^k$, we have
\begin{align}\label{E: k-factorization in function}
  \tksvar{k}{\lambda\cup R_i} &= 
  \tksvar{k}{R_i} \tks{k}{\la}(y|\rot^i(b)).
\end{align}
\end{thm}
 \begin{proof}
We remark that 
$\pi_i(\tksvar{k}{w})=\tks{k}{w}(y|\rot^i(b))$ holds by \eqref{E: pi acts on L by conjugation of finite part} and Remark \ref{R:Wfin action}.
We have \cite[Lemma 2.15]{Tak}
\begin{equation}
x_{\la\cup R_i}=\rot^i(x_\la)x_{R_i}.
\end{equation}
It follows from Proposition \ref{P:factor kappa}
\begin{align*}
\ell_{x_{\la\cup R_i}}
&=\ell_{\kappa_i}\pi_i(\ell_{x_\la})\quad \text{in $\mathbb{L}_{\GL_n}$}.
\end{align*}
From Corollary \ref{C: SL and PGL in functions} (2) and the first remark,
we have \eqref{E: k-factorization in function} in the quotient ring $\hat{\Lambda}_{\PGL_n}$.
We see that \eqref{E: k-factorization in function} holds in $\hLGL$ 
because all the results in \S \ref{SS:k rectangles} hold for $G=\GL_n$
when we use $\tilde{\om}_i^\vee$ (Remark \ref{rem: lift of fundamental weights}) instead of $\om_i^\vee$
(see \cite{IINY} for a more computational approach).
\end{proof}

\begin{ex} 
Let
$k=4$ ($n=5$) and $i=2$. Then $R_2=(2^3).$ Let
$\lambda=(3,1).$ Then $\la\cup R_2=(3,2^3,1).$
$$
  \tksvar{k}{(3,1)\cup R_2} = 
  \tksvar{k}{R_2} \tks{k}{(3,1)}(y|\rot^2(b)).
$$
\end{ex}

\subsection{$SL_n$ and $\PGL_n$ cases}\label{SS: SL and PGL}
In the remaining part of this section we explain how to obtain results for the $\SL_n$ and $\PGL_n$ cases from the $GL_n$ case.

\subsubsection{Root datum for groups of type $A_{n-1}$}
\label{SS:type A root data}

Recall the $G=\GL_n(\C)$ root data from Example \ref{X:GL root data}.

The adjoint type group is 
$\Gad=\PGL_n(\C)=G/\Gm$ and the simply connected group is $\Gsc=\SL_n(\C)\subset G$.
The image 
$\overline{T}$ of $T$ under $G\rightarrow \PGL_n(\C)$ is a maximal torus of $\PGL_n(\C)$, while $\TSL:=T\cap \SL_n(\C)$ is a maximal torus of $\SL_n$.
Let 
$$
a=a_1+\cdots +a_n,\quad
\eps=\eps_1+\cdots +\eps_n
$$
The character and cocharacter groups of $\TSL$ and $\TPGL$ are
given by
\begin{align*}
X(\TSL)&=X(\TGL)/\Z\, a,&
X^\vee(\TSL)&=(\Z\, a)^\perp,\\
X(\TPGL)&=(\Z\,\eps)^\perp,&
X^\vee(\TPGL)&=X(\TGL)/\Z \,\eps,
\end{align*}
where $(\Z\, a)^\perp\subset 
X^\vee(T)$ is the set of all $x\in X^\vee(T)$ such that 
$\langle x,a\rangle=0$, and
$(\Z\, \eps)^\perp\subset 
X(T)$ is the set of all $y\in X(T)$ such that 
$\langle \eps,y\rangle=0$.
We have the following diagram
\[
\begin{tikzcd}
X(\TPGL) = (\Z\eps)^\perp \arrow[d,hookrightarrow] \arrow[r, dotted, leftrightarrow]
& X^\vee(\TSL)=(\Z a)^\perp \arrow[d,hookrightarrow]\\
X(T) = \bigoplus_{i=1}^n \Z a_i \arrow[d,twoheadrightarrow] \arrow[r, dotted, leftrightarrow]& 
X^\vee(T) = \bigoplus_{i=1}^n \Z \varepsilon_i \arrow[d,twoheadrightarrow]\\
X(\TSL) = X(T)/\Z a \arrow[r, dotted, leftrightarrow]& X^\vee(\TPGL) = X^\vee(T)/\Z \varepsilon
\end{tikzcd}
\]
where the horizontal arrows are isomorphisms induced by the noncanonical isomorphism $X(T)\cong X^\vee(T)$.
Let
$
\Phi_\circ, \Phi^\vee_\circ
$ (resp. $\overline{\Phi}, \overline{\Phi}^\vee
$) be the root system and the coroot system of $\SL_n(\C)$ (resp. $\PGL_n(\C)$).
Then $\Phi_\circ$ (resp. $\overline{\Phi}^\vee$ ) is the image of $\Phi$ (resp.  $\Phi^\vee$) under the projection $X(\TGL)\rightarrow X(\TSL)$ (resp.
 $X^\vee(\TGL)\rightarrow X(\TPGL)$).
One sees that ${\Phi}$ (resp. $\Phi^\vee$) is contained in $X(\TPGL)$ (resp. $X^\vee(\TSL)$), and is identified with $\overline{\Phi}$ (resp. $\Phi_\circ^\vee$).
Thus we clearly see that the root datum of $\SL_n(\C)$ and $\PGL_n(\C)$ are dual to each other. 

Set
$$
Q=\Z\Phi, \quad
Q^\vee=\Z\Phi^\vee,\quad
Q_\circ=\Z\Phi_\circ, \quad
Q_\circ^\vee=\Z\Phi_\circ^\vee,\quad
\overline{Q}=\Z\overline{\Phi},\quad
\overline{Q}^\vee=\Z\overline{\Phi}^\vee.
$$
We have
\begin{equation*}
X^\vee(\TSL)=Q_\circ^\vee ,\quad
X(\TPGL)=\overline{Q}.
\end{equation*}
These equalities correspond to the facts that 
$\SL_n(\C)$ is simply-connected, and 
 $\PGL_n(\C)$ is of adjoint type respectively.

\begin{rem}\label{rem: lift of fundamental weights}
The fundamental coweights of $\PGL_n$ are 
\begin{equation*}
\om_i^\vee=\eps_1+\cdots+\eps_i \mod \Z\eps\quad \text{for $1\le i\le n-1.$}
\end{equation*}
They form a $\Z$-basis of $X^\vee(\TPGL)=\bigoplus_{i=1}^n\Z\eps_i/\Z\eps,$ which is 
naturally identified with the coweight lattice of $\GL_n.$
With this identification, the map $\phi$ defined by \eqref{E:Xvee to Pvee} for $G=\GL_n$ is nothing but the natural projection
$$\phi: X^\vee(T)\rightarrow X^\vee(\TPGL).$$
There is a standard lift of $\om_i^\vee\in X^\vee(\TPGL)$
to $X^\vee(\TGL),$ that is, we can define
\begin{equation*}
\tilde{\om}_i^\vee=\eps_1+\cdots+\eps_i \in X^\vee(T)\quad \text{for $1\le i\le n-1$}
\end{equation*}
as a representative for $\om_i^\vee.$
\end{rem}

\subsubsection{Extended affine Weyl groups in type $A$}
Recall that for type $A_{n-1}$, the finite Weyl group $\Wfin$ is the symmetric group $S_n$ and
the affine Wey group is $\hat{W}_{\SL_n}=Q^\vee \rtimes S_n=\langle s_0,s_1,\ldots,s_{n-1}\rangle$ where 
$s_0=t_{\theta^\vee}s_\theta$.
The extended affine Weyl group of $\GL_n$ is
\begin{align*}
\tilde{W}_{\GL_n}&=X^\vee(T) \ltimes S_n=(X^\vee(T)/Q^\vee)\ltimes \hat{W}_{\SL_n}
=\langle \sh, s_0,s_1,\ldots,s_{n-1}\rangle,
\end{align*}
where $X^\vee(T)/Q^\vee \cong \pi_1(\GL_n)\cong \Z$ and $\sh=t_{\eps_1}s_1s_2\cdots s_{n-1}$ is a generator.
We have $\sh s_i =s_{i+1} \sh\;(i\in I_\af=\Z/n\Z).$
The extended affine Weyl group of $\PGL_n$ is
\begin{align*}
\tilde{W}_{\PGL_n}&=X^\vee(\TPGL) \ltimes S_n
=(X^\vee(\TPGL)/\overline{Q}^\vee) \ltimes \hat{W}_{\SL_n}.
\end{align*}
Since $X^\vee(\TPGL)=X^\vee(\TGL)/\Z\eps$, 
$\tilde{W}_{\PGL_n}$ is naturally  a quotient of $\tilde{W}_{\GL_n}$.
In fact, 
$\FG{\PGL_n}=X^\vee(\TPGL)/\overline{Q}^\vee
\cong \Z/n\Z$, and the image of $\sh$ has order $n$ and 
generates $\FG{\PGL_n}.$

\subsubsection{$K_*^{\TSL}(\Gr_{\SL_n})$ and $K_*^{\TPGL}(\Gr_{\PGL_n})$}
\label{SS:SL and PGL}
There are ring homomorphisms $R(T)\rightarrow R(\TSL)$ and $R(\TPGL)\rightarrow R(\TGL)$ that are surjective and injective respectively.

\begin{prop}\label{P: K homology of SL and PGL}
For $w\in \tilde{W}_G^0$, we denote by $\O_{w}^{G}$ the corresponding Schubert class in $\Gr_{G}.$ Then we have the following.
\begin{enumerate}
\item Let $K_*^T(\Gr_{\GL_n})^\circ:=\bigoplus_{x\in \hat{W}_{\SL_n}^0}R(\TGL)\O_{w}^{\GL_n}$.  There is an $R(\TSL)$-Hopf algebra isomorphism 
\begin{align*}
K_*^{\TSL}(\Gr_{\SL_n})&\cong R(\TSL)\otimes_\Rep K_*^{\TGL}(\Gr_{\GL_n})^\circ, \\
 \O_{w}^{\SL_n}&\mapsto 1 \otimes \O_{w}^{\GL_n}\qquad\text{for $w\in \hat{W}_{\SL_n}^0$.}
\end{align*}

\item There is an $R(T)$-Hopf algebra isomorphism 
\begin{align*}
K_*^{T}(\Gr_{\GL_n})/(\O_{\sh^n}-1)&\cong 
\Rep\otimes_{R(\TPGL)}K_*^{\TPGL}(\Gr_{\PGL_n}), \\
\O_{\tau^i v}^{\GL_n}&\mapsto 1\otimes \O_{\pi_i v}^{\PGL_n}\qquad\text{for $i\in\Z$, $v\in \hat{W}_{\SL_n}^0$.}
\end{align*}
\end{enumerate}
\end{prop}
\begin{proof} 
$\Gr_{\SL_n}$ is the identity component of $\Gr_{\GL_n}$. Its torus-fixed points are indexed by $\hat{W}_{\SL_n}^0$. This implies (1). 
We have (2) from Proposition \ref{P: Gr Gad}. 
\end{proof}

Let $\hLSL$
be the $R(\TSL)$-span of 
$\tilde{g}_{w}^{(k)}(y|b)\;(w\in \hat{W}_{\SL_n}^0).$

\begin{cor}\label{cor:Omega/Omega
}
$\Omega(b_i|y)/\Omega(b_j|y)
\in \hLSL.$
\end{cor}
\begin{proof}
From Proposition \ref{prop:Om-Pieri}, $\Omega(b_i|y)/\Omega(b_n|y)
\in \hLSL$ for $1\le i\le k.$ Because $\hLSL$ is invariant under the action of $S_n$, $\Omega(b_i|y)/\Omega(b_j|y)
\in \hLSL.$
\end{proof}

\begin{cor}\label{C: SL and PGL in functions}
We have the following.
\begin{enumerate}
\item  
There is an $R(\TSL)$-Hopf algebra isomorphism 
\begin{align*}
K_*^{\TSL}(\Gr_{\SL_n})&\cong \hLSL ,\\
\O_{w}^{\SL_n}&\mapsto \tilde{g}_{w}^{(k)}(y|b)\quad\text{for $w\in \tilde{W}_{\SL_n}^0$.}
\end{align*}
\item  
There is an $R(T)$-Hopf algebra isomorphism
\begin{align*}
\Rep\otimes_{R(\TPGL)}K_*^{\TPGL}(\Gr_{\PGL_n})&\cong 
\hat{\Lambda}_{(n)}/({\Omega(b_1|y)\cdots \Omega(b_n|y)}-1),\\
1\otimes_{R(\TPGL)}\O_{w}^{\PGL_n}&\mapsto \tilde{g}_{w}^{(k)}(y|b)\quad\text{for $w\in \tilde{W}_{\PGL_n}^0$.}
\end{align*}
\end{enumerate}
\end{cor}
\section{Centralizer family}\label{S:cent}
In this section we realize the equivariant $K$-homology rings 
of type $A$ groups $G$ as the coordinate rings of 
centralizer family $\mathscr{Z}_{G^\vee}$ associated to $G^\vee$.
\subsection{Centralizer family for $\GL_n$}
Consider the matrix $A$ \eqref{E:matrix A}
with entries in $R(\TGL)$.
For $t\in \TGL$ we denote 
by $A(t)$ the element of $\GL_n(\C)$ obtained from $A$ by
evaluating at $t\in \TGL.$   
$B$ denote the Borel subgroup of $\GL_n(\C)$ of upper triangular matrices in $\GL_n(\C)$. 
Define a closed subscheme  $\mathscr{Z}_{\GL_n}$ of  
$\TGL\times B$ by 
 \[
\mathscr{Z}_{\GL_n}= \{(t,Z)\in \TGL\times B\mid ZA(t)=A(t)Z\}.
 \]
 By the 1st projection, $\mathscr{Z}_{\GL_n}$ is a scheme over $\TGL$, which we call the \emph{centralizer family
 for} $\GL_n(\C).$
 \begin{rem}\label{rem:Z}
 It is easy to see that any $n\times n$ matrix $X$ that commutes with $A$
is upper triangular.
So if $X\in \GL_n(\C)$, we have $X\in B.$
 \end{rem}
Let 
$z_{ij}\;(1\le i\le j\le n)$ be the coordinate functions of $B\subset \GL_n(\C)$.  
The commutativity of $Z=(z_{ij})\in B$ and $A$ is equivalent to the relations: 
\begin{equation}
\label{E: central}
(b_i-b_j)z_{ij}=z_{i,j-1}-z_{i+1,j}\quad\text{for $1\le i< j\le n$.}
\end{equation}
Let $R$ be the $\Rep$-algebra
\begin{align*}
   R&:= \Rep[z_{ij}\mid1\le i\le j\le n]/I,
\end{align*}
where $I$ is the ideal generated by the relations \eqref{E: central}.

For an affine scheme $X$ we denote
by $\O(X)$ its coordinate ring.
We have $\O(T)= \Rep$.
Then $\O(\mathscr{Z}_{\GL_n})$ is the localization of $R$ by the multiplicative set generated by $z_{ii}\;(1\le i\le n).$

\begin{lem}\label{lem:flat}
$\mathscr{O}(\mathscr{Z}_{\GL_n})$ is a flat $R(\TGL)$-module. \end{lem}
\begin{proof}
In the ring $R$, we can express each $z_{ij}$
as an $\RepGL$-linear combination of $z_{11},z_{12},\ldots,z_{1n}$ in the following form:
\[
z_{ij}=z_{1,j-i+1}+\sum_{j-i+1<k\le n}c_{ij}^{k}z_{1k},\quad
c_{ij}^k\in \RepGL.
\]
Let \[
f_{ij}=
z_{ij}-(z_{1,j-i+1}+\sum_{j-i+1<k\le n}c_{ij}^{k}z_{1k}).
\]
It is straightforward to show that $I=\langle f_{ij}\;(1\le i\le j\le n)\rangle.$ 
By eliminating $z_{ij}\;(2\le i\le n)$, one sees that $R$ is  
the polynomial ring 
$\RepGL[z_{11},z_{12},\ldots,z_{1n}]$.

Therefore $\O(\mathscr{Z}_{\GL_n})$, given as a localization of $R$, is regular and in particular Cohen-Macaulay. $\Rep$ is clearly regular. 
We apply a criterion of flatness by fibers.
Let $\mathfrak{m}$ be a maximal ideal of $\Rep.$
The fiber ring $\O(\mathscr{Z}_{\GL_n})/\mathfrak{m}\O(\mathscr{Z}_{\GL_n})$ is a
localization of the polynomial ring 
$(\Rep/\mathfrak{m})[z_{11},z_{12},\ldots,z_{1n}]$, and hence has dimension $n.$
Therefore $\O(\mathscr{Z}_{\GL_n})$
is flat over $\Rep$ by 
\cite[Lemma 10.128.1]{Stack}. 
\end{proof}

\begin{prop}\label{prop:W preserve O}
There is an action of $\WGL=S_n\ltimes \bigoplus_{i=1}^n\Z \eps_i$ as ring automorphisms on $\O(\mathscr{Z}_{\GL_n})$.
$S_n$ acts naturally on $\Rep$. For $1\le i\le n-1$, define the action of $s_i$ by
\begin{align}
s_i(z_{jl})&=z_{jl}\; \text{unless $l=i$ or $j=i+1$},\label{eq:siz1}\\
s_i(z_{ii})&=z_{i+1,i+1} \; \text{and $s_i(z_{i+1,i+1})=z_{ii},$} \label{eq:siz2}
\\
s_i(z_{ji})&=z_{ji}+(b_{i+1}-b_{i})z_{j,i+1}\quad \text{for $j<i$} ,\label{eq:siz3}\\
s_i(z_{i+1,l})&=z_{i+1,l}-(b_{i+1}-b_{i})z_{il}\quad \text{for $i<l-1$}.\label{eq:siz4}
\end{align}
The translation elements $t_{\eps_i}$ acts by multiplication
by $z_{ii}$:
\begin{equation}
t_{\eps_i}(f)
=z_{ii}f\quad \text{for $f\in \O(\mathscr{Z}_{\GL_n})\quad( 1\le i\le  n).$}
\label{eq:translation on Z}
\end{equation}
In particular, the action of $s_0$ is given by 
\begin{equation}
s_0(f)= {z_{11}}{z_{nn}}^{-1}s_{\theta}(f)
\quad \text{for $f\in \O(\mathscr{Z}_{\GL_n}).$}
\label{eq: s0 on Z}
\end{equation}
\end{prop}
\begin{proof} 
For $1\le i\le n-1,$
let 
${u}_i=1+(b_{i}-b_{i+1})E_{i+1,i}.$
One sees that 
$u_iAu_i^{-1}$
is equal to $s_i(A)$ by the action of $S_n$ on the entries in $\Rep.$
The matrix $s_i(Z):={u}_iZ{u}_i^{-1}$ clearly 
commutes with $s_i(A).$
Therefore $s_i(Z)$ is an element of $B$ due to Remark \ref{rem:Z} (1). Thus, the entries 
of $s_i(Z)$ satisfies
equations \eqref{E: central} with $b_i$ and $b_{i+1}$ exchanged. This shows 
$A\mapsto s_i(A), Z\mapsto s_i(Z)$
gives an automorphism of $\mathscr{Z}_{\GL_n}.$
It is not difficult to check \eqref{eq:siz1},\eqref{eq:siz2},\eqref{eq:siz3},\eqref{eq:siz4}, and $s_i\;(1\le i \le n-1)$ gives an action of $S_n$.

We define the action of $t_{\eps_i}$ by
\eqref{eq:translation on Z}.
It is straightforward to check 
$w t_{\eps_i}=t_{w(\eps_i)} w$ for $w\in S_n$
and $1\le i\le n.$
\eqref{eq: s0 on Z} holds 
because 
$s_0=t_{\theta^\vee}s_\theta$.
\end{proof}
\begin{rem}
    The same construction on
    the centralizer family in the homology case was given in 
    \cite[Lecture 17]{Pet}.
\end{rem}
\begin{prop}
There is an action of $\KGL$ on $\O(\mathscr{Z}_{\GL_n}).$ 
\end{prop}
\begin{proof}
We first 
show there is an action of $T_i =c(\alpha_i)^{-1}(s_i-1)$ on $\O(\TGL\times B):=\Rep[z_{ij}\,|\,1\le i\le j\le n][z_{ii}^{-1}\,|\, 1\le i\le n].$
Well-definedness is obvious for $1\le i\le n-1$. 
For $f\in \O(\TGL\times B)$,
\[
c(\al_0)z_{nn}T_0(f)=
z_{11}s_\theta(f)-z_{nn}f=
(s_\theta-1)\left(
z_{nn}f
\right).
\]
Because 
$(s_\theta-1)
(z_{nn}f)$ is divisible by $b_1-b_n=c(\al_0)(1-b_n),$
and $1-b_n$ and $z_{nn}$ are invertible, we have $T_0(f)\in \O(\TGL\times B).$

Let $J\subset \O(\TGL\times B)$ be the defining ideal of $\cZ_{\GL_n}.$
By Proposition \ref{prop:W preserve O}, we have seen that 
the action of $\WGL$ on  $\O(\TGL\times B)$ preserves the ideal $J$.
The left $\RepGL$-action clearly preserves $J$.
It remains to show the operators $T_i$ preserve the ideal $J$. By Lemma \ref{lem:flat}, 
$\O(\cZ_{\GL_n})$ is flat over $\RepGL.$ Therefore in particular the map from 
 $\O(\cZ_{\GL_n})\to \O(\cZ_{\GL_n})$
given by the multiplication
with $c(\al_i)$ is injective (since $\Rep\overset{c(\al_i)}{\longrightarrow} \Rep$ is). 
Suppose $f\in J$. If it were the case $T_i(f)\notin J$, then also $c(\al_i)T_i(f)\notin J.$
This is contradiction because 
\[
c(\al_i)T_i(f)=(s_i- 1)f\in J.
\]
Therefore $T_i(f)\in J$.
\end{proof}

 \subsection{
 Parametrizing 
 $\mathscr{Z}_{\GL_n}$ by symmetric functions}
We will construct 
a realization of 
$\O(\cZ_{\GL_n})$ by using symmetric functions.
Let $y=(y_1,y_2,\ldots)$ be the variables as before. 
 Define for $m\ge 0$ and a set of parameters $\{t_1,\ldots,t_i\}$
\begin{equation}
\hat{h}_m(y|t_1,\ldots,t_i):=
\sum_{l=0}^\infty
 h_{l}(t_1,\ldots,t_i)h_{m+l}(y).
\end{equation}

\begin{lem}\label{lem:zij} We have 
\begin{equation}
  \hat{h}_{j-i}(y|b_i,\ldots,b_j)
=e^{a_i+\cdots+a_{j-1}}
\Omega(b_i|y)\cdot \ks{k}{\rho_{j-i}}(y|\rot^i(b))  
\label{eq:hhat}
\end{equation}
for $1\le i\le j\le n.$
\end{lem}
\begin{proof}
It suffices to show that
 \begin{align}\label{E:hhat unshifted}
\hat{h}_{l}(y|b_1,\ldots,b_{l+1})
=e^{a_1+\cdots+a_{l}}
\Omega(b_1|y)\cdot \ks{k}{\rho_{l}}(y|\rot(b))
\qquad\text{for $1\le l\le k$}
\end{align} 
since \eqref{eq:hhat} is obtained from \eqref{E:hhat unshifted}
by setting $l=j-i$ and applying $\rot^{i-1}$.

For $1\le l\le k$, one can verify that
\begin{equation}
T_{l}\cdots T_1(\Omega(b_1|y))=e^{-a_1-\cdots-a_{l}}\hat{h}_l(y|b_1,\ldots,b_{l+1}).
\end{equation}
Thus it remains to show
\begin{equation}
\label{eq:TTOmega}
T_{l}\cdots T_1(\Omega(b_1|y))
=\Omega(b_1|y)\cdot \ks{k}{\rho_{l}}(y|\rot(b)),
\end{equation}
which is shown by by induction on $l.$ If $l=1$
\begin{align*}
T_1(\Omega(b_1|y))
&=\frac{\Omega(b_2|y)-\Omega(b_1|y)}{b_2\ominus b_1}\\
&=\Omega(b_1|y)\left(
\frac{\Omega(b_2|y)/\Omega(b_1|y)-1}{b_2\ominus b_1}
\right) \\
&=\Omega(b_1|y)\rot
\left(
\frac{\Omega(b_1|y)/\Omega(b_n|y)-1}{b_1\ominus b_n}
\right) \\
&=\Omega(b_1|y)\ks{k}{\rho_{1}}(y|\rot(b)).
\end{align*}
For $l\ge 2$, we have
\begin{align*}
T_l\cdots T_1(\Omega(b_1|y))
&=
T_l\left(
\Omega(b_1|y)\ks{k}{\rho_{l-1}}(y|\rot(b))
\right)\\
&=(T_l \Omega(b_1|y))\ks{k}{\rho_{l-1}}(y|\rot(b))
+s_l(\Omega(b_1|y))
T_l(\ks{k}{\rho_{l-1}}(y|\rot(b)))
\\
&=\Omega(b_1|y)\ks{k}{\rho_{l}}(y|\rot(b)),
\end{align*}
because 
 $T_{l}(\ks{k}{\rho_{l-1}}(y|\rot(b)))=\ks{k}{\rho_{l}}(y|\rot(b))$,
 and $\Omega(b_1|y)$ is fixed by $s_l$ for $l\ge 2.$
\end{proof}

\begin{thm}\label{thm:O(Z) and Lambda}
There is a $R(\TGL)$-algebra and $\KGL$-module isomorphism 
 \begin{align*}
\O(\cZ_{\GL_n}) &\isomap \hLGL, \\
{z}_{ij}&\mapsto \hat{h}_{j-i}(y|b_i,\ldots,b_j)\quad\text{for $1\le i\le j\le n$.}
\end{align*}
In particular, $z_{ii}$ is sent to $\Omega(b_i|y)$ for $1\le i\le n$.
The map is also an isomorphism of $\mathbb{K}_{\GL_n}$-modules. 
\end{thm}
\begin{proof}
From the simple identity
\[
h_l(b_i,\ldots,b_{j-1})-h_l(b_{i+1},\ldots,b_{j})
=
(b_i-b_j)
h_{l-1}(b_i,\ldots,b_j),
\]
one sees that 
the functions 
$\hat{h}_{j-i}(y|b_i,\ldots,b_j)$
satisfy the defining equations of $\mathscr{Z}_{\GL_n}$.
Therefore there is a morphism of  $\RepGL$-algebras
$\map_{\GL_n}: \O(\cZ_{\GL_n})\to 
\hL{\RepGL}.$
By Lemma \ref{lem:zij}, the image of $z_{ij}$ is a product of an element of $\Rep\subset\hLGL$, $\Omega(b_i|y)$ which is in $\hLGL$ by Remark \ref{R:Omegas are GL classes}, and $g^{(k)}_{\rho_{j-i}}(y|\omega^i(b)) = \omega^i g^{(k)}_{\rho_{j-i}}(y|b)$ which is in $\hLGL$ since $\hLGL$ is $S_n$-stable.

Let us show the injectivity. 
Let 
$\O(\mathscr{Z}_{\GL_n})^{\Delta}:=R(\TGL)^\Delta\otimes \O(\mathscr{Z}_{\GL_n})$.
Thanks to Lemma \ref{lem:flat}, it suffices to show that the localized homomorphism 
$  \map_{\GL_n}^\Delta:
\O(\cZ_{\GL_n})^\Delta\rightarrow 
\hL{\Rep^\reg}$ is injective.
As an $R(\TGL)^\Delta$-algebra, $\O(\mathscr{Z}_{\GL_n})^{\Delta}$ is generated by
$z_{ii}^{\pm 1} (1\le i\le n)$. Since the image 
of $\{z_{ii}\mid 1\le i\le n\}$ under $\map_{\GL_n}$ is $\{\Omega(b_i|y)\mid 1\le i\le n\}$, which is 
algebraically independent (Lemma \ref{L:Jac}), it follows that $\map_{\GL_n}$ is injectivite.

By construction, the map
$\map_{\GL_n}: \O(\cZ_{\GL})\to \hL{\RepGL}$ is
clearly a $\KGL$-module homomorphism. 
Hence we have 
$\beta_{\GL_n}(T_w (1))=T_w(\beta_{\GL_n}(1))=T_w(1)=\ksvar{k}{w}$
for $w\in \WGLgr.$ This shows in particular that $\map_{\GL_n}$ is surjective.
\end{proof}

\begin{cor}\label{thm:generating} The $\RepGL$-algebra
$K_*^\TGL(\Gr_{\GL_n})$ is generated by 
 $${\rot^j}(\O_{\rho_i})\quad\text{for $1\le i\le n-1,\; 0\le j\le n-1$},\quad
 \O_{\sh^i}^{\pm 1}\quad\text{for $1\le i\le n$}.
 $$
\end{cor}
\begin{proof} As an 
$\RepGL$-algebra
$\O(\cZ_{\GL_n})$ is generated by
$z_{ij}\;(1\le i\le j\le n).$ 
So 
$\hLGL$ is generated by the corresponding elements 
$e^{a_i+\cdots+a_{j-1}}
 \Omega(b_i|y)\cdot \ks{k}{\rho_{j-i}}(y|\rot^i(b)),$  and therefore 
$K_*^\TGL(\Gr_{\GL_n})$ 
is generated by 
$e^{a_i+\cdots+a_{j-1}}
 \rot^{i-1}(\O_{\sh})\cdot {\rot^i}(\cI_{\rho_{j-i}}).$
Since $\cI_{\rho_{j}}=\O_{\rho_j}-\O_{\rho_{j-1}}$,
the corollary holds. 
\end{proof}

\subsection{Centralizer family for $\PGL_n$}
Recall that $\TSL$ is the subgroup of $\TGL$ consisting of the elements of 
determinant $1$.
There is an action of $\mathbb{G}_m$ on $\mathscr{Z}_{\GL_n}$ by scalar multiplication on the matrix $Z.$
The centralizer family of $\PGL_n$ is defined as
$$\cZ_{\PGL_n}:=
(\mathscr{Z}_{\GL_n}|_{\TSL})/\mathbb{G}_m,$$
which is a scheme 
over $\TSL.$

Since the relation \eqref{E: central} is homogeneous in $z_{ij}$, 
$\O(\cZ_{\GL_n})$ is naturally a $\Z$-graded ring such that $z_{ij}$ has degree $1.$
For $m\in \Z$, any element in the degree $m$ part $\O(\mathscr{Z}_{\GL_n})_m$ is written for some $d\ge 0$ as 
$$
\frac{f}{(z_{11}\cdots z_{nn})^d},
$$
where $f$ is represented by a homogeneous polynomial in $z_{ij}$ of degree $nd+m$ with coefficients in $\Rep.$
Then $\O(\cZ_{\GL_n}|_{\TSL})=R(\TSL)\otimes_{\Rep}\O(\cZ_{\GL_n})$ is also $\Z$-graded.
We have
\begin{equation}
\O(\cZ_{\PGL_n})=\O(\cZ_{\GL_n}|_{\TSL})_{0}
\end{equation}
where $\O(\cZ_{\GL_n}|_{\TSL})_{0}$ is the degree $0$ part of $\O(\cZ_{\GL_n}|_{\TSL})$.

\begin{prop}$\O(\cZ_{\PGL_n})$ has a structure of an $R(\TSL)$-Hopf algebra and 
a left $\mathbb{K}_{\SL_n}$-module.
\end{prop}
\begin{proof}
$\cZ_{\PGL_n}$ is a group scheme over $\TSL,$ so it has a natural structure of 
$R(\TSL)$-Hopf algebra.
$R(\TSL)\otimes_\Rep \KGL$ acts on 
$\O(\cZ_{\GL_n}|_{\TSL})$. By a natural map
$\mathbb{K}_{\SL_n}\rightarrow R(\TSL)\otimes_\Rep \KGL$,
$\O(\cZ_{\GL_n}|_{\TSL})$ is a $\mathbb{K}_{\SL_n}$-module.
$\O(\cZ_{\GL_n}|_{\TSL})_0$ is stable by the action of $\mathbb{K}_{\SL_n}$ because 
the action of $\Z[Q^\vee]$ preserves the degree. 
\end{proof}

\begin{thm}\label{thm:beta} There is an isomorphism of $R(\TSL)$-Hopf algebras
$$\O(\cZ_{\PGL_n})\isomap \hLSL$$
such that
\begin{align}
{z_{ij}}/{z_{11}}&\mapsto(-1)^{j-i}\frac{\Omega(b_i|y)}{\Omega(b_1|y)}
g_{i-j}^{(k)}(y|\rot^i(b)).\label{eq:beta}
\end{align}
The map is also an isomorphsm of $\mathbb{K}_{\SL_n}$-modules. 
\end{thm}
\begin{proof}

  As an $R(\TSL)$-algebra, 
$\O(\cZ_{\PGL_n})$ is generated by $z_{ij}/z_{11},$ $1\le i\le j\le n,$
and $z_{11}^n/\det(Z).$
By Lemma \ref{lem:zij} and Theorem \ref{thm:O(Z) and Lambda}, 
we can define a morpshim of $R(\TSL)$-algebra
\[
\beta:
\O(\cZ_{\PGL_n})
\rightarrow
\hL{R(\TSL)}
\]
by the formula
\eqref{eq:beta}. It is clear that 
$\beta$ is $\mathbb{K}_{\SL_n}$-linear. 
Since the ratio $\Omega(b_i|y)/\Omega(b_j|y)$ 
is contained in $\hLSL$ (Proposition  \ref{prop:Om-Pieri}), the image of $\beta$ is contained in $\hLSL$.
Because $\beta$ is $\mathbb{K}_{\SL_n}$-linear, the image of $\beta$ contains 
$\tksvar{k}{w}$ for $w\in \hat{W}_{\SL_n}^0$. Hence 
the image of $\beta$ is $\hLSL$. Clearly $\beta$ is also an $R(\TSL)$-Hopf morphism.
\end{proof}

 \begin{cor}$\Rep$-algebra 
$\hLSL$ is generated by \begin{align*}
    \tks{k}{\rho_i}(y|\rot^j(b))&
\quad \text{for $1\le i\le n-1,\;0\le j\le n-1,$}\\
 \Omega(b_i|y)/\Omega(b_i|y)&\quad \text{for $1\le i,j\le n$}.
\end{align*}
 \end{cor}
 
\subsection{Centralizer family for $\SL_n$}

Let $\cZ_{\SL_n}$ be the closed subscheme of $\cZ_{\GL_n}$ 
consisting of $(A,Z)\in \cZ_{\GL_n}$ such that $\det(Z)=1$. Then we have
\begin{equation}
\O(\cZ_{\SL_n})=\O(\cZ_{GL_n})/\langle z_{11}\cdots z_{nn}-1\rangle.
\end{equation}
\begin{prop} There are isomorphisms of $\Rep$-Hopf algebras
$$R(T)\otimes_{R(\TPGL)}K_{*}^{\TPGL}(\Gr_{\PGL_n})\cong\O(\cZ_{\SL_n})
\cong \hLGL/\langle \Omega(b_1|y)\cdots\Omega(b_n|y)-1\rangle.
$$
\end{prop}
\begin{proof} Since $\O_{\sh^n}$ corresponds to $z_{11}\cdots z_{nn}$
 via the isomorphism $K_*^T(\Gr_{\GL_n})\cong \O(\cZ_{\GL_n})$, 
the first isomorphism is a consequence of Proposition \ref{P: K homology of SL and PGL} (2). One sees that 
the second isomorphism holds by Remark \ref{R: O sh n}.
\end{proof}

\appendix

\section{Infinite rank  analogues: $K$-Molev functions}\label{S:Molev}
In \cite{LLS:backstableK}, the 
$K$-\emph{Molev functions}, indexed by arbitrary partitions $\la\in \Par$, were defined as
representives of the ideal shaves of boundaries of  
Schubert varieties in an infinite dimensional Grassmannian.
These functions are $R(T_\Z)$ deformations of the dual stable Grothendieck polynomials $g_\la(y)$ of Lam and Pylyavskyy \cite{LP},
where $R(T_\Z)=\Z[e^{\pm a_i}\mid i\in \Z]$.
The functions presented here are related to the ones in \cite{LLS:backstableK}
by the substitution of variables specified in \cite[\S 2.8]{LLS:backstableK}.

Let $S_\Z$ be the group of permutations $w$ of $\Z$ such that $w(i)=i$ for all but finitely many $i\in \Z.$
$S_\Z$ is presented by the generators $s_i\;(i\in \Z)$ and the relations $s_i^2=1$ and $s_is_{i+1}s_i=s_{i+1}s_is_{i+1}$ for $i\in \Z$
and $s_is_j=s_js_i $ for $|i-j|\ge 2.$
Let $S_\Z^0$ be the set of \emph{$0$-Grassmann elements}, the set of $w\in S_\Z$ such that $w(i)<w(i+1)$ for all $i\in \Z\setminus\{0\}.$
There is a bijection $\Par\cong S_\Z^0$ denoted by $\la\mapsto w_\la$ \cite[\S 2.1]{LLS:backstable}.
\begin{ex}
For $\la=(4,2,2,1)$ we have $w_\la=s_{-3}(s_{-1}s_{-2})(s_0s_{-1})(s_2s_1s_0).$
\end{ex}

$S_\Z$ acts on $R(T_\Z)=\Z[e^{\pm a_i}\mid i\in \Z]$ such that $s_i$ exchanges $e^{a_i}$ and $e^{a_{i+1}}$.
We set $b_i=1-e^{-a_i}$ for $i\in \Z.$
We extend the action of $s_i$ for $i\in\Z$ to $\hL{R(T_\Z)}$. 
Define the action of $s_0$ 
on $\hL{R(T_\Z)}$ by
$$
(s_0 f)(y|b)=\frac{\Omega(b_1|y)}{\Omega(b_0|y)}f(y|s_0(b)),
$$
and for $i\ne 0$
$$
(s_if)(y|b)=f(y|s_i(b)).
$$
Let $
\alpha_i=a_i-a_{i+1}
$ for $i\in \Z.$
Define the endomorphisms of $\hL{R(T_\Z)}$ 
\begin{equation}
\label{E: TiDi infty}
    \hat{T}_i={c(\alpha_i)}^{-1}{(s_{i}-1}),\quad
\hat{D}_i=\hat{T}_i+1.
\end{equation}

For $\la\in \Par$, define
\begin{equation}
   {g}_\la(y|b):= \hat{T}_{w_\lambda} (1),\quad
    \tilde{g}_\la(y|b):= \hat{D}_{w_\lambda} (1).
\end{equation}
$g_\la(y|b)$ is a double version of the dual stable Grothendieck polynomial of Lam and Pylyavsky \cite{LP}. 

Define the ring homomorphism $R(T_\Z)\rightarrow \Rep$ sending $b_i$ to $b_{i\,\mathrm{mod}\, n}$.
We denote the image of $f(y|b)\in \hL{R(T_\Z)}$ 
under this map by 
$f(y|b^{(n)}).$
\begin{rem}\label{R: k-small}
Under the induced map 
$\hL{R(T_\Z)}\rightarrow \hL{R(T)}$, for
$0\le i\le n-1$, the actions of $\hat{D}_i,\hat{T}_i$ descend to
the actions of $D_i,T_i\in \KGL$ on $\hL{\Rep}.$  
\end{rem}

A partition is \emph{$k$-small} if it is contained in a 
$j \times (k+1-j)$-rectangle for some $1\le j \le k$.

\begin{lem}\label{L :k-small}
Suppose $\la\in \Par^k$ is $k$-small. 
Then $x_\la\in \WSLgr$ is obtained from $w_\la$ by 
replacing $s_i$ with $s_{i\,\mathrm{mod}n}$.
\end{lem}
\begin{proof} 
Let $\le_L$ denote the left weak order on $\WSLgr$.
We claim that  for any $k$-small $\la,\mu\in \Par^k$, 
 $x_\la\le_L x_\mu$ is equivalent to $\la\subset \mu$.
Let $\la\in \Par^k$ be $k$-small. Then for any the box in $\la$ the hook length is less than  $n$. 
This implies that the $n$-core of $\la$ equals $\la.$ Then the claim follows from 
\cite[Proposition 25]{LM:JCTA}.
\end{proof}

\begin{prop}
If $\la\in\Par^k$ is $k$-small, then 
\begin{equation}
\tksvar{k}{\la}=\tilde{g}_\la(y|b^{(n)}),\quad
\ksvar{k}{\la}={g}_\la(y|b^{(n)}).
\end{equation}
\end{prop}
\begin{proof}
This is clear from Remark \ref{R: k-small} and Lemma \ref{L :k-small}.
\end{proof}

\section{Double $k$-Schur functions}
\label{S:double k schur}
The equivariant homology of $\Gr_{\SL_n}$
and the double $k$-Schur functions were studied in \cite{LS:double Schur}.
The functions used in \cite{LS:double Schur} are based on the most common definition of
factorial Schur function in the literature. Unfortunately this formulation does not
naturally agree with equivariant localization values.
For the equivariant (co)homology of the infinite Grassmannian this situation was remedied in \cite{LLS:backstable}, where the geometric factorial Schurs and corresponding version of Molev's dual Schurs were described. 
The purpose of this section is to give a version of double $k$-Schur function which is 
consistent with localization and also with the double $K$-theoretic $k$-Schur functions of this paper and the backstable double Schubert polynomials of \cite{LLS:backstable}.

\subsection{Double $k$-Schur functions}
We define double $k$-Schur functions indexed by $w\in \WGLgr$.
Let $G=\GL_n(\C)$ and $S=\Z[a_1,\ldots,a_n]=H_T^*(pt)$.
The group $\WGL$ (see Examples \ref{X:GL translation elements} and
\ref{X: extended affine Weyl type A}) acts on $\hL{S}$ as follows.
$S_n$ acts by permuting $a_1,\dotsc,a_n$.
The element $t_{\pm\epsilon_i}$ acts by multiplication by $\Omega(a_i|y)^{\pm1}$ where
\begin{equation}
\Omega(a_i|y)=\prod_{j=1}^\infty(1-a_i y_j)^{-1}\qquad\text{for $1\le i\le n$.}
\end{equation}
For $1\le i\le n-1$ we define the operator $A_i = (a_i-a_{i+1})^{-1}(1-s_i)$
acting on the $a_i$ variables in $\hLa_S$. We define the operator $A_0 = (a_n-a_1)^{-1}(1 - t_{\theta^\vee} s_\theta)$ where $\theta^\vee=\epsilon_1-\epsilon_n$ and $s_\theta$ exchanges $1$ and $n$. For $0\le i\le n-1$ the $A_i$ operators satisfy the braid relations for the affine symmetric group.
For $w\in \WGLgr$ the double $k$-Schur function is defined by
\begin{equation}
s_w^{(k)}(y|a)=A_{w}(1)\in \hL{S}.
\end{equation}

\begin{prop} There is an $S$-Hopf algebra isomorphism
$\bigoplus_{w\in \WGLgr} S s_w^{(k)}(y|a) \cong H_*^T(\Gr_{\GL_n})$ such that 
$s_w^{(k)}(y|a)$ corresponds to 
the Schubert class $\sigma_w\in  H_*^T(\Gr_{\GL_n})$.
\end{prop}

For $\la\in \Par^k,$
we denote $s_{x_{\la}}^{(k)}(y|a)$ by
$s_\la^{(k)}(y|a)$. 
Let $s_\la^{(k)}(y)$
denote the $k$-Schur function \cite{LLMSSZ}. 
Let $\Lambda_d$ denote the degree $d$ part of $\Lambda.$
\begin{prop}
Let $\la\in \Par^k.$ Then $s_\la^{(k)}(y|a)$ has the form
$$
s_\la^{(k)}(y|a)
=s_\la^{(k)}(y)+\sum_{i\ge 1}c_i(a)f_{|\la|+i}(y),
$$
where $c_i(a)\in S$ is homogeneous of degree $i$, and 
$f_{|\la|+i}(y)\in \Lambda_{|\la|+i}.$
\end{prop}

\begin{rem}
Let 
$\hat{S}:=\mathbb{Q}[[a_1,\ldots,a_n]]$. Define a ring homomorphism 
$\Rep\rightarrow \hat{S}$ by $e^{a_i}\mapsto \sum_{j=0}^\infty \frac{a_i^j}{j!}.$
Let $\mathfrak{m}=(a_1,\ldots,a_n)$ the maximal ideal of $\hat{S}.$
We regard $g_w^{(k)}(y|b)$
as element of $\hL{\hat{S}}$. Then we have
\begin{equation}
g_w^{(k)}(y|b)=
s_w^{(k)}(y|a)+\sum_{i\ge 0}\tilde{c}_i(a){f}_{\ell(w)+i}(y)
\end{equation}
where $\tilde{c}_i(a)\in \mathfrak{m}^{i+1},\; {f}_{\ell(w)+i}(y)\in \Lambda_{\ell(w)+i}.$ 
\end{rem}
\begin{rem}
$s_\la^{(k)}(y|a)$
satisfies the $k$-rectangular factorization. 
\end{rem}

\subsection{Molev dual Schur functions}
The version in \cite{LLS:backstable}  of Molev's dual Schur functions \cite{Molev:comult}
are an infinite rank analogue 
of the double $k$-Schur functions $s_w^{(k)}(y|a)$. We consider
$\Z[a_i\mid i\in \Z]$ with the natural action of $S_\Z$
with $s_i$ exchanging $a_i$ and $a_{i+1}$ for all $i\in\Z$.
Let us denote the action by $s_i^a.$
For $i\in\Z$ define
\begin{align*}
\hat{A}_i &=(a_i-a_{i+1})^{-1}(1-s_i^a), \\
A_i &= \Omega(a_0|y)^{-1} \hat{A}_i \Omega(a_0|y). 
\end{align*}
For $i\ne0$ $A_i = \hat{A}_i$ and
\begin{align*}
  A_0 &= (a_0-a_1)^{-1}(1 - \Omega(a_1|y) \Omega(a_0|y)^{-1} s_0^a).
\end{align*}
The $\hat{A_i}$ and therefore the $A_i$, satisfy the braid relations.
For $w\in S_\Z$, $\hat{A}_w$ is defined
by using a reduced expression of $w$.
For $w\in S_{\Z}^0$, the Molev dual Schur function is defined by \cite[\S 8]{LLS:backstable}
\begin{equation}
\hat{s}_w(y|a)=A_{w}(1).
\end{equation}
For $\la\in \Par$, we denote $\hat{s}_{w_\la}(y|a)$ by $\hat{s}_\la(y|a).$
Define the ring homomorphism $\Z[a_i\mid i\in \Z]\rightarrow S$ sending $a_i$ to $a_{i\,\mathrm{mod}\, n}$.
We denote the image of $f(y|a)\in \hL{\Z[a_i\mid i\in \Z]}$ 
under this map by 
$f(y|a^{(n)}).$
\begin{prop}
If $\la\in\Par^k$ is $k$-small, then 
\begin{equation}
s^{(k)}_\la(y|a)=
\hat{s}_{\la}(y|a^{(n)}).
\end{equation}
\end{prop}

We will prove an analogue of Molev's Jacobi-Trudi formula
\cite[Proposition 3.9 and Corollary 3.10]{Molev:comult}.
 For $j\ge 1$ ( cf. \eqref{E:single row}, \eqref{E:single column}) define $\rho_j,\rho_j'\in\WGL$ by
 \begin{align*}
     \rho_j=s_{j-1}\cdots s_1 s_0,\quad
     \rho_j'=s_{-j+1}\cdots s_{-1}s_0.
 \end{align*}
 Let the operator $\sh$ act on $\hL{S_\Z}$ by sending $a_i$ to $a_{i+1}$ for all $i\in \Z.$
 \begin{lem}\label{L:symmetry in a0 and a1}
    For $i,j\ge 1$ such that $i\ne j$, $\hat{s}_{\rho_i}(y|\sh^{-j}a)$ is invariant under the action of $s_0^a$; that is, it is symmetric in $a_0$ and $a_1$.
 \end{lem}
 \begin{proof}
 It is easy to see that $\hat{s}_{\rho_i}(y|a)$ is symmetric in $a_1,a_2,\ldots, a_i$. Therefore  
for $i,j\ge 1$ such that $j<i$, 
$\hat{s}_{\rho_i}(y|\sh^{-j}a)$ is symmetric in $a_0$ and $a_1$.
If $i,j\ge 1$ satisfies $j>i$, $\hat{s}_{\rho_i}(y|\sh^{-j}a)$ does not depend on $a_0$ and $a_1$, so 
it is symmetric.
\end{proof}

\begin{prop}
Let $\lambda$ be a partiton and
$\la'$ be the conjugate of $\la$.
 Then
\begin{align}
\hat{s}_{\la}(y|a)&=
\label{eq:JT1}
\det(\hat{s}_{\rho_{\la_i+j-i}}(y|{\sh^{-j+1}}a))_{i,j}    \\
&=
\label{eq:JT2}
\det(\hat{s}_{\rho_{\la_i'+j-i}'}(y|\sh^{j-1}a))_{i,j}.
\end{align}
\end{prop}
\begin{proof}
Let $\la\in \Par$, and $i\in \Z$. We define that $\la$ is $i$-addable if 
$s_i w_\la>w_\la$ and $s_iw_\la\in S^0_\Z.$ Then
the corresponding partition is denoted by $\la+i.$
Let $M_\la=(\hat{s}_{\rho_{\la_i+j-i}}(y|{\sh^{-j+1}}(a)))_{i,j}$.
Suppose $\la$ is $i$-addable. We will show 
$$
 A_i \det M_\la=
\det M_{\la+i} 
$$
From this identity, we have $\det M_\la=A_{w_\la}(1)$ so \eqref{eq:JT1} is obtained. 
We only show this for the case of $i=0$; the other cases are easier.

For $j\ge 1$, we will show
\begin{equation}\label{eq: A0 acts in Pieri}
A_0 (\hat{s}_{\rho_{j-1}}(y|\sh^{-j+1}a))=\hat{s}_{\rho_{j}}(y|\sh^{-j+1}a).
\end{equation}
Let $\sh$ acts on $\hL{S}$ by sending $f(y|a)$ to $\Omega(a_1|y)f(y|\sh a).$ Then 
$\sh$ is invertible and $\sh A_i =A_{i+1}\sh.$
We have (cf. \eqref{eq:sh acts on K-k-Schur})
\begin{equation}
\label{eq: sh hat s}
\sh^j (\hat{s}_{\la}(y|a))=\sh^j(1)\cdot \hat{s}_\la(y|\sh^ja).
\end{equation}
If $j=1$ \eqref{eq: A0 acts in Pieri} holds by definition. Suppose $j\ge 2$, then 
$$
A_0 \sh^{-j+1}A_{j-2}\cdots A_1A_0=\sh^{-j+1}A_{j-1} A_{j-2}\cdots A_1A_0
$$
together with \eqref{eq: sh hat s} 
implies \eqref{eq: A0 acts in Pieri}.

Suppose $\la=(\la_1,\ldots,\la_l)$ is $0$-addable. We may assume $l\ge 2.$ 
Then 
for some $s\ge 2$, we have $\la_s=s-1$ and $\la_{s-1}\ge s$.
We see from Lemma \ref{L:symmetry in a0 and a1} that all the entries of $M_\la$ except for the ones in the $s$-th row are $s_0^a$-invariant.
We expand the determinant along the $s$-th row
$$
\det M_\la=\sum_{j=1}^{\ell(\la)}(-1)^j\hat{s}_{\rho_{j-1}}(y|\sh^{-j+1}a) \det M_\la^{(j)},
$$
where $\det M_\la^{(j)}$ is the $(s,j)$-minor. 
Note that for $f,g\in \hL{S}$ such that $f$ is $s_0^a$-invariant, we have  $A_0(fg)=f A_0(g)$. Therefore we have
\begin{align*}
A_0 \det M_\la&=\sum_{j=1}^{\ell(\la)}(-1)^j\left(A_0\hat{s}_{\rho_{j-1}}(y|\sh^{-j+1}a)\right) \det M_\la^{(j)} \\
&=\sum_{j=1}^{\ell(\la)}(-1)^j\hat{s}_{\rho_{j}}(y|\sh^{-j+1}a) \det M_\la^{(j)}\quad \text{by \eqref{eq: A0 acts in Pieri}}\\
&=\det M_{\la+i}.
\end{align*}
\eqref{eq:JT2} is obtained from \eqref{eq:JT1} by applying an automorphism 
of $\hL{S}$ such that $\hat{s}_{\la}(y|\sh^i a)\mapsto \hat{s}_{\lambda'}(y|\sh^{-i}a)$
(cf. \S \ref{sec: k-conj}).

\end{proof}
\begin{rem} The operator $\sh$ on $\hL{S}$ defined in the proof generates a group $L$ isomorphic to $\Z.$
We can consider the semidirect $L\rtimes S_\Z$ such that $\sh s_i =s_{i+1} \sh$ for $i\in\Z$,
and introduce an infinite rank version of the nil-Hecke algebra. 
See \cite{LLS:backstable} for more details on such construction. 
\end{rem}

\end{document}